\documentclass[a4paper, 12pt]{article}

\usepackage[sort&compress]{natbib}
\bibpunct{(}{)}{;}{a}{}{,} 

\usepackage{amsthm, amsmath, amssymb, mathrsfs, multirow, url, subfigure}
\usepackage{graphicx} 
\usepackage{ifthen} 
\usepackage{amsfonts}
\usepackage[usenames]{color}
\usepackage{fullpage}
\usepackage{tikz}

\theoremstyle{plain} 
\newtheorem{thm}{Theorem}
\newtheorem{cor}{Corollary}

\theoremstyle{definition}

\theoremstyle{remark} 
\newtheorem{ex}{Example}

\newcommand{\prob}{\mathsf{P}}

\newcommand{\unif}{{\sf Unif}}
\newcommand{\nm}{{\sf N}}

\newcommand{\gam}{{\sf Gamma}}

\newcommand{\chisq}{{\sf ChiSq}}

\newcommand{\RR}{\mathbb{R}}

\newcommand{\TT}{\mathbb{T}}
\newcommand{\T}{\mathcal{T}}

\newcommand{\eps}{\varepsilon}

\newcommand{\prior}{\mathsf{Q}}

\newcommand{\cred}{\mathscr{C}}

\newcommand{\lPi}{\rotatebox[origin=c]{180}{$\mathsf{\Pi}$}} 
\newcommand{\uPi}{\mathsf{\Pi}}

\title{Probabilistic approximations of possibility measures and fiducial inference for discrete models}
\author{Ryan Martin\footnote{Department of Statistics, Purdue University, {\tt martinrg@purdue.edu}} \quad and \quad Sahil S. Patel\footnote{Department of Statistics, North Carolina State University, {\tt sspate27@ncsu.edu}}}
\date{\today}

\begin{document}

\maketitle 

\begin{abstract}
Discrete-data models are fundamental and seemingly simple, yet have lurking complexities that create challenges when it comes to reliable, probabilistic uncertainty quantification.  Indeed, Fisher famously failed to offer a version of his fiducial argument for discrete-data problems.  A suggested alternative to Fisher's fiducial argument is to construct a provably reliable possibilistic inferential model (IM) and then suitably approximate it by a probability.  In discrete-data models, however, the aforementioned complexities manifest themselves as jump discontinuities in the IM's contour, which is an obstacle to obtaining a probabilistic approximation.  This paper provides a characterization of this probabilistic approximation and, importantly, shows that it generally contains data-dependent atoms, which no Bayesian or fiducial solutions can have.  We also prove that this probabilistic approximation has an exact, finite-sample reliability guarantee and is asymptotically efficient.

\smallskip

\end{abstract}

\section{Introduction}
\label{S:intro}

The {\em inferential model} (IM) framework, introduced by \citet{imbasics, imbook} and more recently summarized in \citet{imreview}, offers fully conditional, Bayesian-like uncertainty quantification with strong finite-sample frequentist-like reliability guarantees.  
What brings the best-of-both-worlds within reach is the fact that, while Bayesian and fiducial uncertainty quantification is probabilistic, the IM framework's is {\em possibilistic} \citep[e.g.,][]{dubois.prade.book, dubois2015.where, dubois2006}.  In fact, possibility theory is arguably the ``right'' brand of uncertainty quantification for statistical inference, as there is no loss of generality or efficiency in viewing it through a possibility-theoretic lens; see Section~\ref{S:background}. 

Despite the promising theory, a practical challenge is that the Monte Carlo methods advanced over the last few decades \citep[e.g.,][]{liang.liu.carroll.book, casella, jsliu2001} are specifically designed for probabilistic uncertainty quantification---that is, for evaluating probabilities and expected values---and, therefore, are not immediately applicable for IM's brand of possibilistic uncertainty quantification.  So, new methods and/or insights are needed.  Breakthrough progress was made in \citet{immc}, where it was shown that the IM's possibilistic output could be suitably approximated by a probability measure, to facilitate Bayesian-like sampling on the parameter space, and then transformed back to a possibility measure for inference with little or no material loss in the IM's reliability guarantees.  Roughly, a mixture characterization was established for those probability measures inside the possibility measure's credal set and, in turn, this leads to a characterization of the so-called {\em inner probabilistic approximations} as the ``most diffuse'' member of said credal set.  We review these developments in Section~\ref{SS:imreview}.  

While the advancements just described are completely general, the practical suggestions based on this theory may be challenging to implement accurately when the IM's output is not ``smooth'' in a sense to be defined below.  For example, discrete data problems such as binomial inference---one of the first statistical problems given careful treatment \citep[e.g.,][]{laplace} and called the ``fundamental problem of practical statistics'' by \citet{pearson1920}---are ones that lack the aforementioned smoothness.  These are obviously very important, practical applications, so a motivation behind the new developments in the present paper is to fill some foundational and computational gaps when it comes to handling discrete-data models.  Below we summarize the paper's key contributions.  

First, in Section~\ref{S:poss.to.prob}, we prove a general characterization of the inner probabilistic approximation of a possibility measure whose contour has jump discontinuities.  This characterization has details in common with other cases that covered in the literature (see Appendix~\ref{A:inner}), but the result in Theorem~\ref{thm:jumpy} is new to our knowledge.  An important and distinguishing detail is that, when the contour has jump discontinuities, the corresponding inner probabilistic approximation will typically have atoms, i.e., single points in the support that are assigned positive probability.  

Second, Section~\ref{S:discrete} proves that, for discrete data models, the possibilistic IM's contour has jump discontinuities.  Then Theorem~\ref{thm:jumpy} implies that inner probabilistic approximations of the possibilistic IM's output must include (data-dependent) atoms.  This is counterintuitive because the parameter space typically contains an open set, so we intuitively expect absolutely continuous distributions for uncertainty quantification.  

Following the recommendation in \citet{reimagined}, our third contribution is to use the possibilistic IM's inner probabilistic approximation as an alternative to fiducial \citep[e.g.,][]{zabell1992, fisher1973}, generalized fiducial \citep[e.g.,]{hannig.review}, or default-prior Bayes solutions \citep[e.g.,][]{berger.objective.book}.  What is particularly interesting about this discrete data context is that there are no agreed-upon solutions: Fisher famously failed to provide a fiducial argument for discrete data models, and other related efforts only offer large-sample frequentist performance guarantees.  One explanation for this shortfall is that there is no realistic way in which these solutions can have the apparently necessary data-driven atoms.  But the proposed probabilistic approximation of the possibilistic IM has atoms and, as we show in Section~\ref{S:fiducial}, in addition to having certain finite-sample calibration properties, is also asymptotically optimal in the sense that it satisfies the same Bernstein--von Mises theorem that the Bayesian and fiducial solutions satisfy.  

The paper concludes in Section~\ref{S:discuss} with a brief summary and discussion of some related extensions and open problems.  Additional details, including proofs of the theorems and other cases not covered in the main text, are provided in the Appendix.

\section{Background}
\label{S:background}

\subsection{Possibility theory}
\label{SS:possibility}

Let $\TT$ denote a generic space, e.g., a Euclidean space $\TT = \RR^d$ for $d \geq 1$.  While not immediately necessary for the possibilistic formulation, when we speak of probability distributions defined on $\TT$, we will assume that $\TT$ is equipped with a $\sigma$-algebra of subsets, such as the Borel $\sigma$-algebra, on which it is supported.  

An important feature that probability and possibility theory have in common is that both are fully determined by scalar-valued function.  In probability theory, that function is a density or mass function, and probabilities of events are defined via integration.  In possibility theory, that function is a {\em contour} and possibilities are defined via optimization.  More formally, let $\pi: \TT \to [0,1]$ be a function.  If it satisfies the normalization condition, $\sup_{t \in \TT} \pi(t) = 1$, analogous to the more familiar summation/integral normalization for mass/density functions in probability theory, then it is a (possibility) contour.  The corresponding {\em possibility measure} $\uPi$ is defined via optimization as 
\begin{equation}
\label{eq:maxitive}
\uPi(H) = \sup_{t \in H} \pi(t), \quad H \subseteq \TT.
\end{equation}
The associated {\em necessity measure} is $\lPi(H) = 1 - \uPi(H^c)$, but we will focus on $\pi$ and $\uPi$.  In applications, $\uPi(H)$ is interpreted as a measure of the possibility or plausibility of the hypothesis $H$ in light of the available information (encoded in $\pi$).
Importantly, only small values of $\uPi(H)$ carry inferential weight: knowing $H$ is {\em plausible} contributes little, but knowing $H$ is {\em implausible} means we have learned something.  

As suggested by the interpretation above, the inequality $\lPi(H) \leq \uPi(H)$ holds for all $H$, with strict inequality except for trivial $H \in \{\varnothing, \TT\}$.  As such, we can think of $\lPi$ and $\uPi$ as special cases of a more general kind of lower and upper probability, respectively, like those that can be found in the belief function literature \citep[e.g.,][]{shafer1976, dempster2008, cuzzolin.book} and the imprecise probability literature \citep[e.g.,][]{walley1991, lower.previsions.book, imprecise.prob.book}.  This connection to imprecise probability theory is important, but we emphasize that, while general imprecise probability can be quite complex, possibility theory is no more complicated than probability theory.  

General upper probabilities like $\uPi$ have an associated {\em credal set} that contains all those ordinary/precise probabilities that are dominated by $\uPi$:
\begin{equation}
\label{eq:cred0}
\cred(\uPi) = \{ \prior \in \text{probs}(\TT): \text{$\prior(H) \leq \uPi(H)$ for all measurable $H$} \}. 
\end{equation}
Note that the probabilities in $\cred(\uPi)$ can be discrete, continuous, or mixtures of both.  This reinforces the above point that only small values of $\uPi(H)$ carry inferential weight, since it implies that small probability is assigned to $H$ by each element in the credal set $\cred(\uPi)$.  The simplicity of possibility measures compared to other brands of upper probabilities allows for more concrete characterizations of the credal set.  Indeed, if $\uPi$ is a possibility measure, then \citet{cuoso.etal.2001} proved that 
\begin{equation}
\label{eq:char0}
\prior \in \cred(\uPi) \iff \prior(C_\alpha) \geq 1-\alpha, 
\end{equation}
where $C_\alpha = \{t \in \TT: \pi(t) > \alpha\}$ is the {\em $\alpha$-cut} of $\pi$.  That is, $\prior$ belongs to $\cred(\uPi)$ if and only if $\pi(T)$ is stochastically no smaller than $\unif(0,1)$ as a function of $T \sim \prior$; this clearly has some close connections to p-values, which is relevant later.  Credal sets are also important when it comes to approximating possibilities by probabilities, a focus of the present paper.  Further details about the connections between probabilities and possibilities through the credal set are summarized in Appendix~\ref{AA:background}.

\subsection{Possibilistic IMs}
\label{SS:imreview}

An inferential model (IM) is a model or procedure that defines a mapping from data, statistical model, etc.~to a quantification of uncertainty about relevant unknowns.  The {\em possibilistic} IM in consideration here is, as described in detail by \citet{imreview}, one whose output takes the form of a (data-dependent) possibility measure.  More specifically, suppose we have observable data $X$ and a posited model $\{\prob_\theta: \theta \in \TT\}$ for said data indexed by a parameter space $\TT$.  The working assumption is that $X \sim \prob_\Theta$, where $\Theta$ is the uncertain ``true'' parameter.  Given $X=x$, the goal is to quantify uncertainty about $\Theta$ in a reliable way.  As we will argue below, the possibilistic IM achieves this goal.  Below we give only the details necessary to follow the arguments in the upcoming sections, but some additional details are provided in Appendix~\ref{AA:background2}. 

\begin{figure}[t]
\begin{center}
\subfigure[Exponential data]{\scalebox{0.55}{\includegraphics{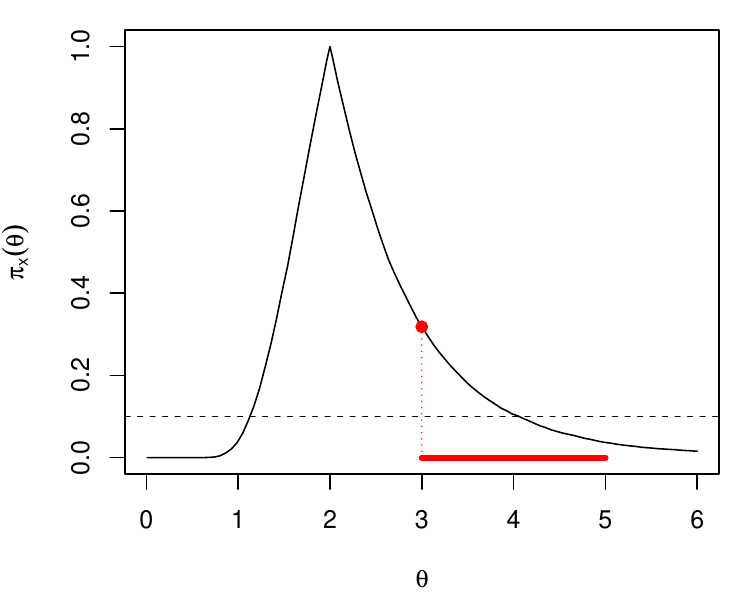}}}
\subfigure[Bernoulli data]{\scalebox{0.55}{\includegraphics{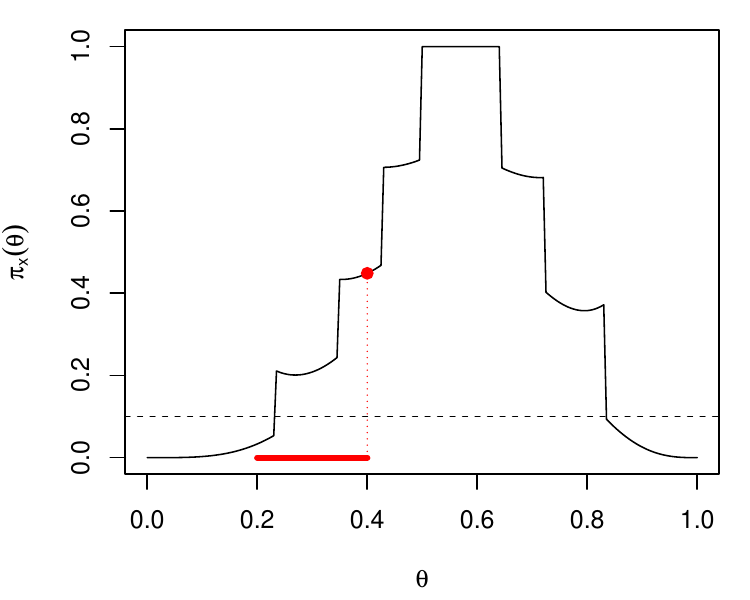}}}
\end{center}
\caption{Plots of the IM contour function for two different models and data sets. Red horizontal line represents a hypothesis $H$ and the red dot corresponds to $\uPi_x(H)$ determined by optimization. The horizontal dashed line at $\alpha=0.1$ determines the 90\% confidence interval derived from the IM output.}
\label{fig:toy}
\end{figure}

Following \citet{imreview}, the possibilistic IM construction may proceed as follows.  If $p_\theta$ is the density or mass function corresponding to the model $\prob_\theta$, then $\theta \mapsto p_\theta(x)$ is the likelihood function based on $X=x$.  Let $\hat\theta_x$ be the/a maximum likelihood estimator, and define the relative likelihood function as
\[ R(x,\theta) = \frac{p_\theta(x)}{\sup_\vartheta p_\vartheta(x)} = \frac{p_\theta(x)}{p_{\hat\theta_x}(x)}. \]
This relative likelihood can be replaced by a different ranking functions, and this is often recommended in situations where, say, the model involves nuisance parameters \citep[e.g.,][]{martin.partial3, imreview}.  But the likelihood-based rankings above are sufficiently general for our purposes here.  Then the possibilistic IM construction proceeds by defining the (data-dependent) possibility contour as 
\begin{equation}
\label{eq:contour}
\pi_x(\theta) = \prob_\theta\{ R(X,\theta) \leq R(x,\theta) \}, \quad \theta \in \TT,
\end{equation}
where ``$X$'' in the above display is a random variable having distribution $\prob_\theta$ and ``$x$'' is the observed data.  That this defines a genuine possibility contour is easy to see: of course, it takes values in $[0,1]$ and it attains the value 1 at $\theta=\hat\theta_x$.  Then the IM's possibility or upper probability is defined via optimization just like as in \eqref{eq:maxitive}:
\[ \uPi_x(H) = \sup_{\theta \in H} \pi_x(\theta), \quad H \subseteq \TT. \]
Figure~\ref{fig:toy} shows two IM contours, both based on iid samples of size $n=7$, from an exponential and Bernoulli distribution, respectively.  In addition to the contour, both show a hypothesis $H$ and the corresponding possibility $\uPi_x(H)$ defined via optimization, and the 90\% confidence interval determined by cutting the contour at $\alpha=0.1$. 

The right-hand side of \eqref{eq:contour} resembles a p-value corresponding to a likelihood ratio test, and this implies strong reliability or {\em strong validity} of the IM output:
\begin{equation}
\label{eq:svalid}
\sup_\theta \prob_\theta\{ \uPi_X(H) \leq \alpha \, \text{for some $H$ with $H \ni \theta$} \} \leq \alpha, \quad \alpha \in [0,1]. 
\end{equation}
In words, this says that assigning a low upper probability to a hypothesis $H$ that is true is a similarly low-probability event---uniformly in hypotheses $H$.  Since we are likely to dismiss or reject those $H$ for which $\uPi_x(H)$ is small, this property ensures that such errors can be controlled in a frequentist sense.  Moreover, \eqref{eq:svalid} implies that the contour's level sets, $C_\alpha(x) = \{\theta: \pi_x(\theta) \geq \alpha\}$, are frequentist confidence sets in the sense that 
\begin{equation}
\label{eq:coverage}
\sup_\theta \prob_\theta\{ C_\alpha(X) \not\ni \theta \} \leq \alpha, \quad \alpha \in [0,1]. 
\end{equation}
It is important to emphasize, however, that this goes beyond p-values, significance testing, etc.  First, the IM provides fully conditional, possibilistic uncertainty quantification without prior distributions---default or otherwise---with strong reliability guarantees.  Further details about the connection to Bayesian, frequentist, and fiducial inference can be found in \citet{imreview, reimagined, imchar, martin.isipta2023} and Section~\ref{S:fiducial} below.  Second, $\uPi_x$ can be extended beyond just assigning possibilities to hypotheses: for instance, \citet{imdec.ext} use Choquet integration with respect to $\uPi_x$ for decision-making under uncertainty.  Finally, under mild conditions \citep[see][App.~G]{imreview}, for any valid frequentist procedure, e.g., any confidence set procedure that satisfies \eqref{eq:coverage}, there exists a possibilistic IM whose contour satisfies \eqref{eq:svalid} and has $\alpha$-cuts that agree with the given confidence set.


One simple but not-practically-useful way in which the reliability described above can be achieved simply by being overly conservative: if the contour is diffuse, then $\uPi_X(\cdot)$ tends to be large, so there are no non-trivial practical hypotheses that can be refuted; similarly, a diffuse contour means wide confidence sets, so no non-trivial parameter values can be ruled out.  But the above IM construct does not achieve its validity properties by being overly conservative.  Indeed, the {\em possibilistic Bernstein--von Mises theorem} of \citet{imbvm.ext} establishes that the above IM solution is asymptotically optimal in the same sense that maximum likelihood estimators and Bayesian posteriors are, under the same regularity conditions.  Section~\ref{S:fiducial} below provides further details on this.

%
%
%

\section{Probabilistic approximations of possibilities}
\label{S:poss.to.prob}

\subsection{Setup and overview}

Inspired by \citet{dubois.etal.2004} and others, the aim of the present section (and the companion Appendix~\ref{A:inner}) is to provide explicit descriptions of a ``best'' approximation of a possibility by a probability.  Such approximations are useful for both theory and computation,. Since the properties established here hold for all possibility measures, we adopt the generic notation from Section~\ref{SS:possibility}, $\pi$ and $\uPi$, omitting the data subscript.  

We start by introducing the desired approximation.  Following \citet{immc, reimagined}, the jumping off point is the characterization \eqref{eq:char0}: each member $\prior$ of the credal set $\cred(\uPi)$ assigns probability at least $1-\alpha$ to the $\alpha$-cuts $C_\alpha = \{t: \pi(t) > \alpha\}$.  So, in order for a probability $\prior^\star$ to represent both the tails and center of $\uPi$, it makes sense to consider the ``most diffuse'' element of $\cred(\uPi)$, i.e., the one that assigns probability {\em exactly} $1-\alpha$ to $C_\alpha$.  This definition is difficult to work with, but a more refined notion was offered in \citet{immc, reimagined}, one based on a mixture characterization of the contents of $\cred(\uPi)$.  An equivalent formulation, also mentioned briefly in the above references, turns out to be quite powerful and will be our working definition here.  That is, an {\em inner probabilistic approximation} $\prior^\star$ of $\uPi$ is an element of the credal set $\cred(\uPi)$ such that 
\begin{equation}
\label{eq:char}
\prior^\star\{ \pi(T) \leq \pi(t) \} = \pi(t), \quad t \in \TT.
\end{equation}
This definition captures the aforementioned p-value-like property. Appendix~\ref{AA:background} provides an interpretation of $\prior^\star$ as an inverse of the so-called {\em probability-to-possibility transform}.


Characterizations of $\prior^\star$ have appeared elsewhere in the literature, but we are not aware of a systematic treatment.  The newfound interest and applications in statistics suggests it is worth fleshing out a few relevant general results.  This section focuses exclusively on the case where $\TT = \RR$ and the contour $\pi$ has jump discontinuities, since this is most relevant to the specific case in Section~\ref{S:discrete}, but Appendix~\ref{A:inner} covers two additional cases that are also quite relevant: the case where $\TT$ is countable and the case where $\TT = \RR^d$, with $d \geq 1$, and the contour is smooth.

\subsection{Case: contour has jump discontinuities}
\label{SS:jumps}

For technical reasons discussed below, we focus our attention here on the case where $\TT$ is a subset of $\RR$.  Suppose that the contour $\pi: \TT \to [0,1]$ has the following properties:
\begin{itemize}
\item has a modal interval (could be a singleton), on which $\pi$ is constant equal to 1;
\vspace{-2mm}
\item non-decreasing from the left towards the modal interval, differentiable with the exception of finitely many jump discontinuities; 
\vspace{-2mm}
\item non-increasing to the right away from the modal interval, differentiable with the exception of finitely many jump discontinuities;
\vspace{-2mm}
\item and the upper level sets $\{t: \pi(t) \geq \alpha\}$ are closed.
\end{itemize} 
This more-or-less describes the typical contour obtained by the possibilistic IM solution when data are discrete; see Section~\ref{S:discrete} for more details on this case.

\begin{figure}[t]
\begin{center}
\scalebox{0.55}{\includegraphics{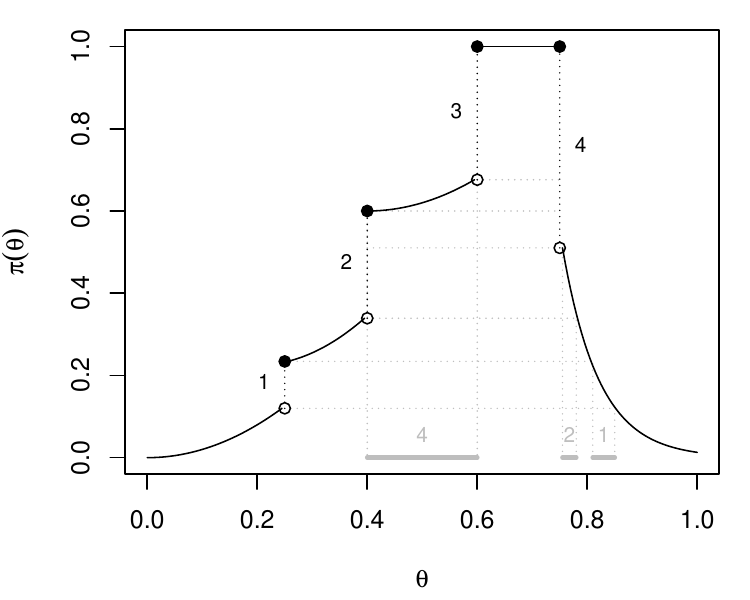}}
\end{center}
\caption{Plot of a generic possibility contour $\pi$ with jump discontinuities.  In this case, there are $J=4$ jump points $(t_1,t_2,t_3,t_4)$ and the vertical jumps in $\pi$ are labeled 1--4 (black).  The shadow intervals $\mathcal{T}_j$ corresponding to jumps $t_j$, for $j \in \{1,2,4\}$, are denoted by the gray line segments on the horizontal axis, labeled accordingly.}
\label{fig:jumpy}
\end{figure}

Roughly speaking, Theorem~\ref{thm:jumpy} below establishes that, if $\uPi$ is a possibility measure with contour that satisfies the above properties, then an inner probabilistic approximation $\prior^\star$ of $\uPi$ is a suitable mixture of discrete and absolutely continuous parts. More formally, the probability distribution, $\prior^\star$, is described via its {\em Lebesgue decomposition}, $\prior^\star = \prior_\text{\sc c}^\star + \prior_\text{\sc d}^\star$, into absolutely continuous and discrete parts, respectively.  Towards this, partition the domain of $\pi$ into finitely many disjoint open intervals on which the function is differentiable, separated by points at which $\pi$ jumps.  If there is a single mode, with $\pi$ increasing before and decreasing after, then split that interval around the mode into two intervals.  Label the jump points as $t_1 < t_2 < \cdots < t_J$.  For the case illustrated in Figure~\ref{fig:jumpy}, there are $J=4$ jump points located at 0.25, 0.4, 0.6, and 0.75, respectively.  

We start with the absolutely continuous part.  Note that jumps on one side of the mode cast {\em shadow intervals} on the other side of the mode on the $\theta$-axis.  In Figure~\ref{fig:jumpy}, two jump points, $t_1$ and $t_2$, on the left of the mode cast shadows on the right, and the one jump on the right of the mode, $t_4$, casts a shadow on the left; as discussed below, jump point $t_3$ casts no shadow (or it casts a singleton shadow).  The idea is that when the contour jumps at, say $t_1$, on the left, there is a range, say $\mathcal{T}_1$, of values on the right---marked by the gray interval labeled ``1'' in Figure~\ref{fig:jumpy}---over which $\theta$ must traverse in order for the change in $\pi(\theta)$ to match that at the jump $t_1$.  That is, $\pi$ evaluated at the left and right endpoints of $\mathcal{T}_1$ equal $\pi(t_1+)$ and $\pi(t_1-)$, respectively.  The jump at $t_3$ has a singleton shadow and hence can be omitted; this cannot affect the absolutely continuous part of $\prior^\star$.  Let $\mathcal{T} = \bigcup_{j=1}^J \mathcal{T}_j$ be the union of the jump shadow intervals.  Then the absolutely continuous part $\prior_\text{\sc c}^\star$ of $\prior^\star$ has a density defined as 
\[ q_\text{\sc c}^\star(\theta) = \begin{cases} |\dot\pi(\theta)| & \theta \in \mathcal{T} \\ \tfrac12|\dot\pi(\theta)| & \theta \not\in \mathcal{T}. \end{cases} \] 
The absolute value accounts for the fact that $\pi$ is non-decreasing to the left of the mode and non-increasing to the right.  The factor $1/2$ accommodates ties, as discussed in Appendix~\ref{SS:common}, and is absent on the shadow intervals because there are no ties.  

Next, the discrete part.  As the reader might guess, the masses that determine the discrete part are concentrated at the jump points $t_1,\ldots,t_J$ of $\pi$.  Counterintuitively, however, not all of the jump points get assigned positive mass. Like in the finite case of Appendix~\ref{SS:finite}, the masses assigned correspond to differences, namely, 
\[ q_\text{\sc d}^\star(t_j) = \Bigl| \pi(t_j) - \sup_{\theta: \pi(\theta) < \pi(t_j)} \pi(\theta) \Bigr|, \quad j=1,\ldots,J. \]
In words, the probability mass assigned to a jump point is the difference between the possibility assigned at that jump point and the maximal possibility among all those that are strictly smaller---note that this difference can be zero. Again, if there are ties, i.e., two distinct jump points $t_i$ and $t_j$ with $\pi(t_i)=\pi(t_j)$, then the smaller of the two masses as defined above is divided between them.  For the case in Figure~\ref{fig:jumpy}, first note that zero probability mass is assigned to $t_1=0.2$, since the supremum of the contour over shadow interval ``1'' is equal to $\pi(0.2)$.  Next, at $t_2=0.4$, the mass assigned is 
\[ q_\text{\sc d}^\star(0.4) = \pi(0.4) - \pi(0.75+) \approx 0.089, \]
the difference between the contour at $t_2$ and the supremum over shadow interval ``2.''  The contour values at $t_3=0.6$ and $t_4=0.75$ are equal, so there is a tie; the mass assigned is 
\[ q_\text{\sc d}^\star(0.6) = q_\text{\sc d}^\star(0.75) = \tfrac12 \{1 - \pi(0.6-)\} \approx 0.162, \]
where the $1/2$ factor divides the mass equally between the two jump points. So, $\prior_\text{\sc d}^\star$ assigns positive mass to three of the four jump points totally, roughly, $0.089 + 2 \times 0.162 = 0.413$; the residual mass, 0.587, is captured by the absolutely continuous part. 

The following theorem establishes that the $\prior^\star$ just described is indeed an inner probabilistic approximation of the given possibility measure $\uPi$ with contour $\pi$.  

\begin{thm}
\label{thm:jumpy}
Let $\uPi$ be a possibility measure with contour function $\pi$ that satisfies the above properties. Then $\prior^\star$ defined above is an inner probabilistic approximation of $\uPi$. 
\end{thm}

We end this section with two remarks.  First, since the mixture characterization of the inner probabilistic approximations offered in \citet{immc} is general, it must be that the more concrete characterization given here is a special case.  The key point is that mixture representation allows for variations in how the allocation of mass on the boundaries of the contour's $\alpha$-cuts, $\alpha \in [0,1]$, and how the marginal distribution for $\alpha$ is spread over $[0,1]$.  So, the concrete representation in Theorem~\ref{thm:jumpy} corresponds to sometimes splitting mass between both endpoints of the $\alpha$-cut and sometimes concentrating all the mass on only one endpoint, depending on whether the contour intersects the horizontal line at level $\alpha$ twice or only once; when the range of the contour is disconnected, so that there are some $\alpha$ values that are never reached by the contour, the gaps are covered by the atoms in $\prior_\text{\sc d}^\star$ and these manifest as atoms in the marginal distribution of $\alpha$. 

Second, while the focus here is possibility measures supported on the real line, the message here surely applies more generally: jumps in the possibility contour will generally imply atoms in the probabilistic approximation.  We do not yet have a satisfactory resolution to the multiparameter case, but see Sections~\ref{S:discrete} and \ref{S:discuss} for further discussion. 


\section{Discrete-data possibilistic IMs}
\label{S:discrete}

\subsection{Probabilistic approximations}
\label{SS:approx}

It is interesting---perhaps even striking---that, for discrete-data models with a continuous, one-dimensional parameter, such as binomial and Poisson, the probabilistic approximation of the possibilistic IM is a distribution with both continuous and discrete components.  That is, the IM's probabilistic approximation has atoms with positive probability mass.  With one caveat, explained below, this follows from Theorem~\ref{thm:jumpy} and the fact that, for discrete-data models, the IM contour has jumps. We have already visually observed these jumps in Figure~\ref{fig:toy}(b), so this is not unexpected.  

The aforementioned caveat is that, at least in some cases, the contour can actually {\em increase} temporarily as $\theta$ moves away from the mode.  See, for example, around the value $\theta=0.25$ in Figure~\ref{fig:toy}(b).  This non-monotonic behavior is undesirable but not easily remedied; \citet{balch.binomial} managed to remove this in the binomial case by replacing the relative likelihood in \eqref{eq:contour} with something else.  This non-monotonicity is undesirable because it could lead to a disconnected confidence ``interval,'' which is awkward, difficult to interpret, etc.  What is often done in such cases is to fill in the missing part to make the disconnected set an interval, thus avoiding the aforementioned complications.  In our context, this amounts to replacing the non-monotonic contour with the smallest monotonic contour that dominates it, which we will do implicitly in what follows without mentioning it or modifying the notation.  With that inconsequential substitution, the conditions of Theorem~\ref{thm:jumpy} are met and we can proceed as described next; for an explanation of why this substitution is ``inconsequential,'' see Appendix~\ref{AA:background2}. 

\begin{thm}
\label{thm:discrete.ipa}
Let $\prob_\theta$ be a discrete-data model indexed by a continuous, scalar parameter $\theta \in \TT \subseteq \RR$.  If the corresponding likelihood function is continuous in $\theta$, then the probabilistic approximation of the possibilistic IM is a distribution with both continuous and discrete components, i.e., it has atoms with positive probability mass.  
\end{thm}

While the focus here is on scalar parameters, the phenomenon is more general.  Indeed, jumps in the possibility contour correspond to points $\theta$ for which there exists different data sets $y$ with $R(y,\theta)=R(x,\theta)$.  A difficulty here is that actually finding the jump points is vastly more difficult in higher-dimensional cases than it is in scalar parameter cases, and we do not yet have a satisfactory solution.  However, in the scalar parameter case, finding the jump points is not too difficult, as we explain next.

\subsection{Computation}

In the scalar-parameter models under consideration here, finding the jumps is relatively easy.  Indeed, one can simply enumerate all the other possible data sets $y$ and solve the equation $R(y,\theta) = R(x,\theta)$ directly.  Remember that ``data sets $y$'' really means just the distinct values of the minimal sufficient statistic; so, in the Bernoulli data case, there are only $n$ equations to solve: one for each element in $\{0,1,\ldots,n\} \setminus \{x\}$.  More formally, suppose the model consists of $n$ iid observations from an exponential family distribution with natural or canonical parameter $\theta$, so that the relative likelihood is of the form 
\[ R(x,\theta) = \exp\bigl[ (\theta - \hat\theta_x) \, s(x) - n\{A(\theta) - A(\hat\theta_x)\} \bigr], \]
where $s$ and $A$ are determined by the model.  Then $R(x,\theta) = R(y,\theta)$ simplifies to 
\[ (\theta - \hat\theta_x) \, s(x) + n A(\hat\theta_x) = (\theta - \hat\theta_y) \, s(y) + n A(\hat\theta_y), \]
which is linear in $\theta$, and, consequently, the jump point solution is 
\[ \theta = \frac{\hat\theta_x \, s(x) - \hat\theta_y \, s(y) - n \{A(\hat\theta_x) - A(\hat\theta_y)\}}{s(x) - s(y)}. \]
Of course, if we can find the jump points for exponential families in canonical form, then we can immediately do the same for different parametrizations.  And beyond exponential families, the solutions can be found numerically.  

Once the jump points have been identified, the solution suggested by Theorem~\ref{thm:jumpy} is within reach.  Between the jumps, we need the contour's derivative, which is rather simple.  That is, between the jump points, the set $\mathcal{R}(x,\theta)$ in the summation \eqref{eq:pi.sum} is unaffected by small changes in $\theta$, so we only need to differentiate the summand:
\begin{equation}
\label{eq:gradient}
\dot\pi_x(\theta) = \sum_{z \in \mathcal{R}(x,\theta)} \Bigl\{ \frac{\partial}{\partial \theta} \log p_\theta(z) \Bigr\} \, p_\theta(z). 
\end{equation}
The logarithmic derivative is just the model's score function and, just like the contour itself, this derivative can be numerically evaluated without errors due to Monte Carlo or approximating derivatives by finite difference ratios. It remains to (a)~decide whether a $1/2$ factor is needed or not, and (b)~evaluate the discrete mass assigned to each jump point.  Both of these steps require care, but the details are tedious so we simply refer the interested reader to \url{https://martinrg07.github.io/Codes/imdisfid_binom.R}.

\begin{figure}[t]
\begin{center}
\subfigure[Discrete part, $n=20$]{\scalebox{0.39}{\includegraphics{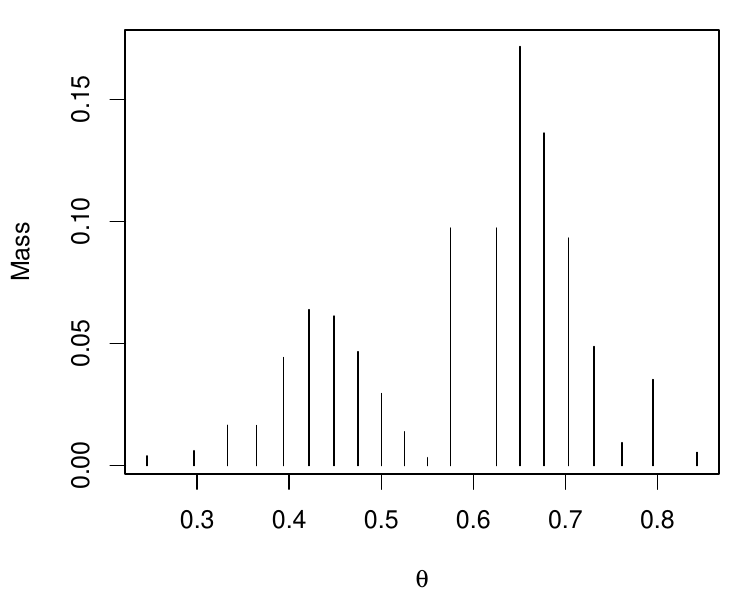}}}
\subfigure[Discrete part, $n=40$]{\scalebox{0.39}{\includegraphics{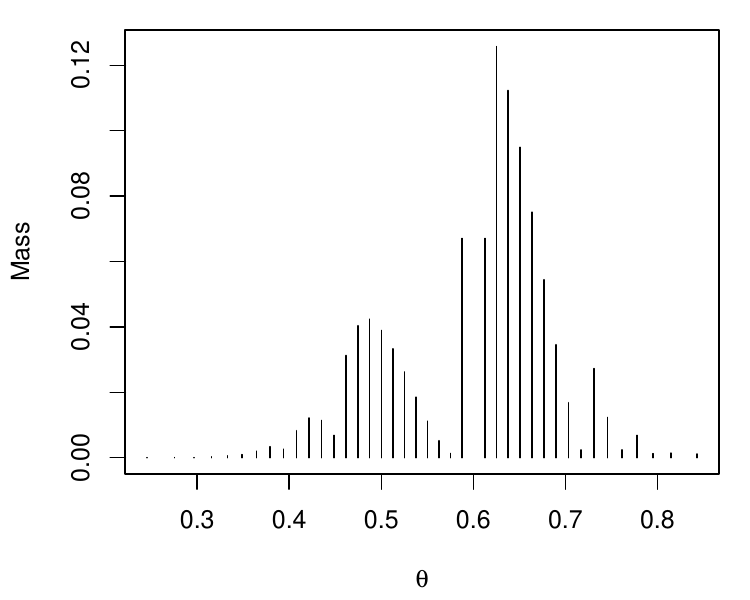}}}
\subfigure[Discrete part, $n=80$]{\scalebox{0.39}{\includegraphics{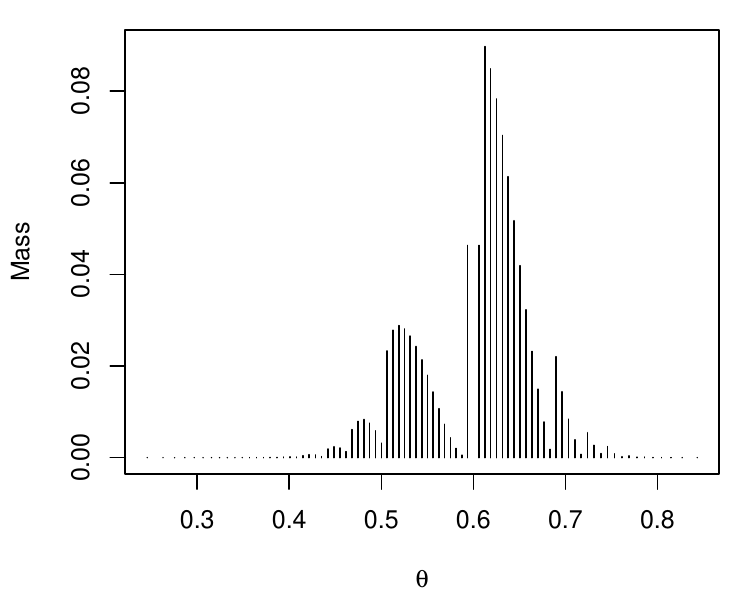}}}
\subfigure[Continuous part, $n=20$]{\scalebox{0.39}{\includegraphics{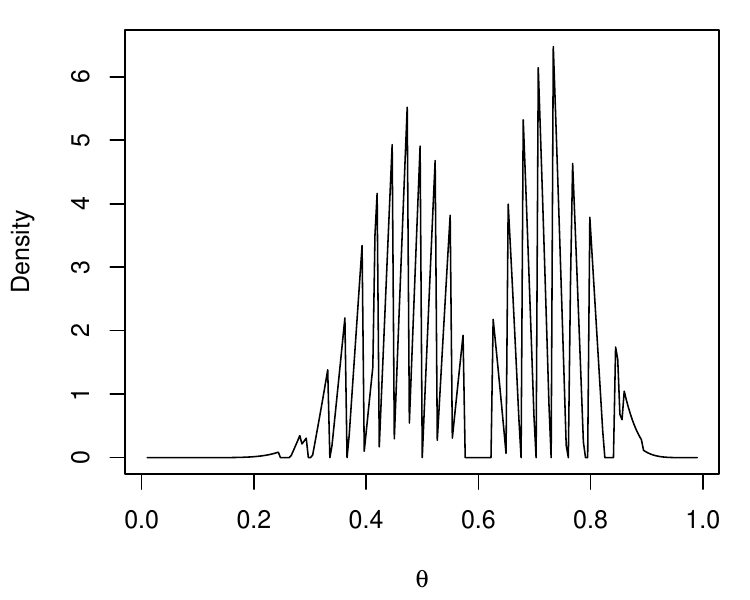}}}
\subfigure[Continuous part, $n=40$]{\scalebox{0.39}{\includegraphics{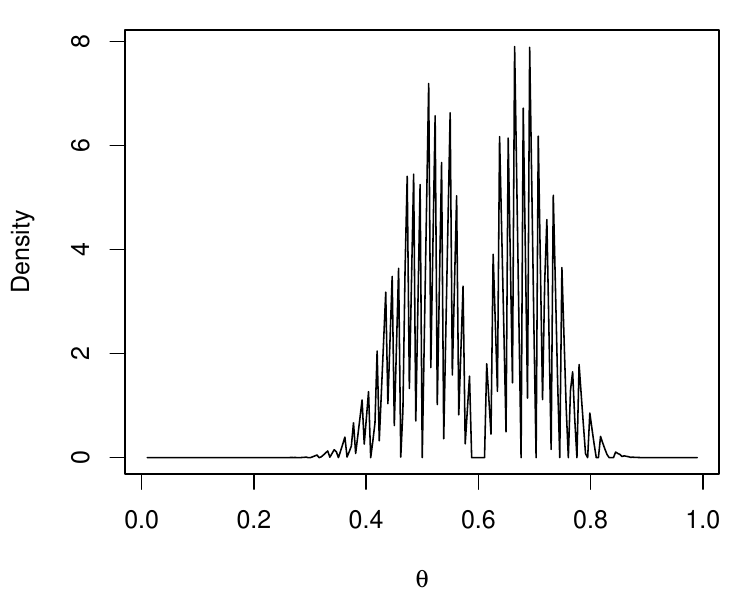}}}
\subfigure[Continuous part, $n=80$]{\scalebox{0.39}{\includegraphics{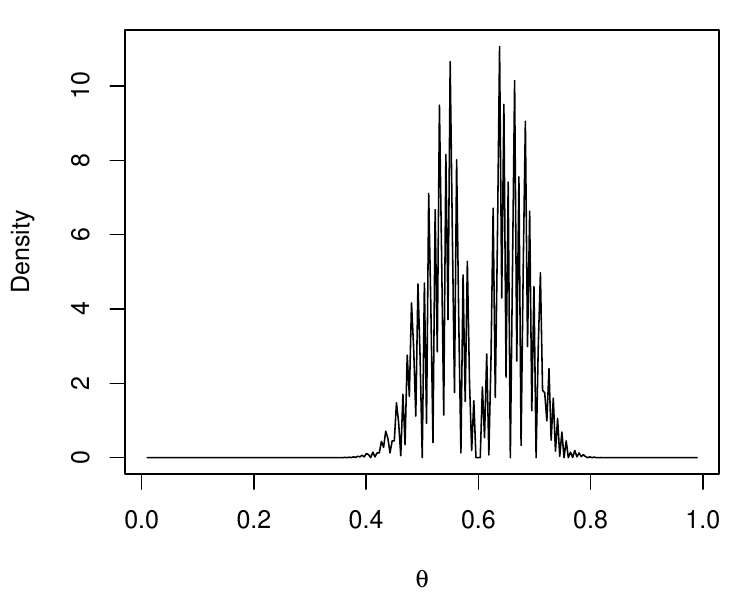}}}
\end{center}
\caption{Plots of the discrete parts (top row) and continuous parts (bottom row) of the inner probabilistic approximation in the binomial model based on three data sets: all three have $\hat\theta = 0.6$, with varying sample sizes $n$.}
\label{fig:discrete}
\end{figure}

As an illustration, consider the binomial model.  Since this is an exponential family, the jump points can be determined in closed-form.  Similarly, between jump points, the the binomial score function is well-known and can be plugged into \eqref{eq:gradient}.  Figure~\ref{fig:discrete} shows both the discrete and continuous parts of the possibilistic IM's inner probabilistic approximation in the binomial model based on three different data sets: all three have the same maximum likelihood estimator, $\hat\theta=0.6$, but with different sample sizes.  One point that is clear from these plots is that the continuous part is very rough, but with an overall shape similar to that of the discrete part.  In particular, both the continuous and discrete parts flatten out in a neighborhood of the maximum likelihood estimator, as expected based on the contour's plateau there.  One point that is not clear from the plot is that the discrete part is overwhelmingly dominant: the total mass assigned to the discrete part is 0.89, 0.93, and 0.97 for $n=20$, 40, and 60, respectively.  So, the apparently very rough continuous part plays at most a very minor role in the approximation, which is perhaps counterintuitive: recall that typical Bayes posterior and fiducial distributions are fully continuous, whereas the probabilistic approximation here is almost entirely discrete.  

\section{The new fiducial argument}
\label{S:fiducial}

Following the proposal in \citet{reimagined}, a new brand of fiducial argument can be based on use of the above-described inner probabilistic approximation of the possibilistic IM.  The proposal is incredibly simple to state: instead of using the posited model to directly construct a data-dependent probability distribution (using the pivot-based argument offered by Fisher, the generalized version of Hannig, etc.),~first find the provably reliable possibilistic IM solution and then return a suitable inner probabilistic approximation thereof.  The intuition behind this proposal is that the IM solution itself has many desirable reliability properties, but its possibility-theoretic framework and structure is unfamiliar; so, if a more familiar probability-theoretic solution is desired, then one should use a probability distribution sufficiently close to the IM solution that it inherits many of the desirable reliability properties.  This proposal is a general one, but the focus in \citet{reimagined} is on continuous-data models: he shows a number of desirable properties of the inner probabilistic approximation, e.g., it agrees with default-prior Bayes in location--scale and other group invariant models, but relevant questions about its use in discrete-data models remain unanswered.  The present section aims to answer these. 

While discrete-data models are among the most fundamental, reliable probabilistic inference remains a challenge.  On the Bayesian side, the models for which there exists a default prior distribution such that the corresponding posterior returns credible sets that are exact confidence sets are effectively only those with group invariance, such as location--scale, which are continuous.  So, in discrete-data cases, Bayesians only get probability matching asymptotically, as sample size increases to infinity \citep[e.g.,][]{welch.peers.1963, datta.ghosh.1995}.  On the fiducial side, as discussed in \citet{zabell1992} and elsewhere, Fisher famously offered no satisfactory fiducial solution for discrete-data problems, though he did eventually advance something, which he hints at in later editions of {\em Statistical Methods and Scientific Inference} \citep[e.g.,][]{fisher1973}.  The solution offered in \citet{fisher.underworld} is cryptic, but our understanding is that he was unable to determine a single probability distribution representing {\em the} fiducial solution, so ambiguity remains.  Similar non-uniqueness issues afflict the generalized fiducial solutions of \citet{hannig.review} and \citet{murph.etal.fiducial} for discrete-data models.  Despite their obvious differences, what these traditional approaches have in common is that they quantify uncertainty about $\Theta$ using a {\em continuous} probability distribution.  Such a choice is quite natural, given that $\Theta$ takes values in an open subset of Euclidean space, but obviously not necessary.  Moreover, no practically justifiable Bayesian solution can have discrete components: the only way the posterior can have discrete components is if the prior itself has discrete components, and there is no {\em a priori} justification for treating some parameter values differently from others.  But Theorem~\ref{thm:discrete.ipa} reveals that reliable probabilistic uncertainty quantification, in the sense of being exactly calibrated to the sampling model and at least asymptotically efficient, apparently must have discrete components, in stark contrast with what all the other attempts at probabilistic uncertainty quantification offer.  The apparently critical importance of discrete components in the uncertainty quantification for discrete-data models explains the difficulties with the traditional methods that are constrained to be continuous, and highlights the benefits of the proposed IM-based approach.  

Of course, the fact that the proposed solution---our probabilistic approximation of the possibilistic IM---has discrete components is not enough to justify its use; the location of and mass assigned to those atoms need to be ``right'' in some sense.  Below we offer two results along these lines.  The first establishes that certain credible sets constructed based on our proposed probabilistic uncertainty quantification are exactly frequentistly calibrated for all sample sizes.  The second shows, via a version of the celebrated Bernstein--von Mises theorem, that provable reliability does not come at the expense of inefficiency. 

\begin{thm}
\label{thm:coverage}
Let $f=f_x$ be a function and consider the corresponding $100(1-\alpha)$\% credible set obtained via $f$-level sets:
\[ C_\alpha^f(x) = \{\theta \in \TT: f_x(\theta) > c_\alpha^f(x)\}, \]
where $c_\alpha^f(x)$ is chosen so that 
\[ \prior_x^\star\{ C_\alpha^f(x)^c \} \equiv \prior_x^\star(\{ \theta: f_x(\theta) \leq c_\alpha^f(x) \}) = \alpha. \] 
If $f=f_x$ is comonotone with $\pi_x$, i.e., $f_x$ is an increasing function of $\pi_x$, then the $100(1-\alpha)$\% credible set $C_\alpha^f(X)$ is also a $100(1-\alpha)$\% confidence set. 
\end{thm}

\begin{figure}[t]
\begin{center}
\subfigure[$n=20$]{\scalebox{0.39}{\includegraphics{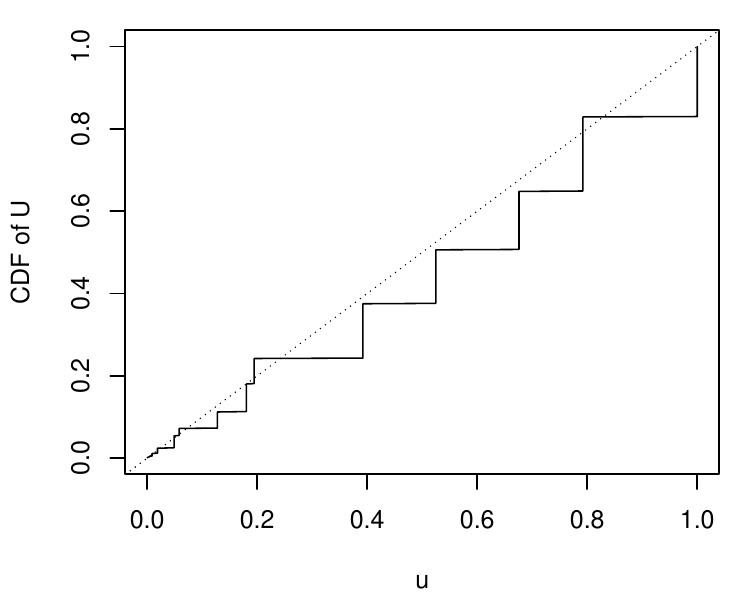}}}
\subfigure[$n=40$]{\scalebox{0.39}{\includegraphics{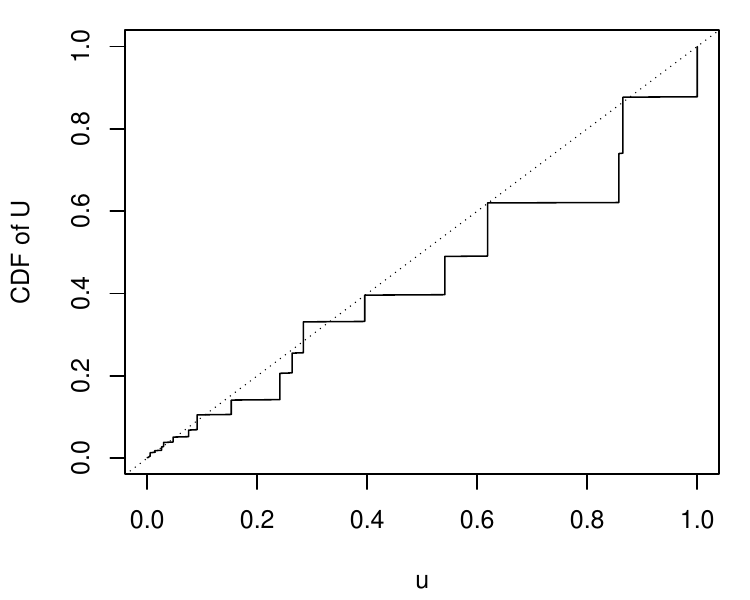}}}
\subfigure[$n=80$]{\scalebox{0.39}{\includegraphics{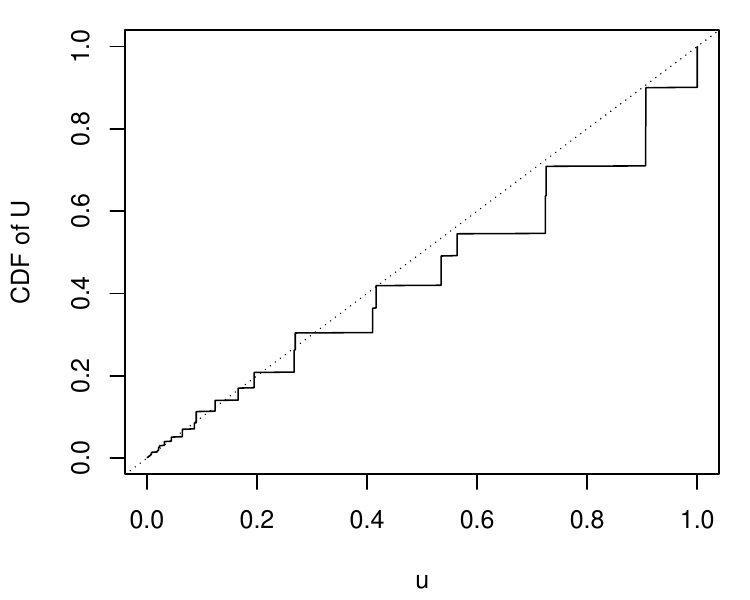}}}
\end{center}
\caption{Plots of the distribution function of $U(X,\Theta)$ based on 1000 simulated data sets from the binomial model with $\Theta=0.6$ for different values of $n$. Diagonal line is the $\unif(0,1)$ reference.}
\label{fig:coverage}
\end{figure}

When $\prior_x^\star$ is absolutely continuous, one might consider taking $f=f_x$ to be the posterior density, which yields the commonly used highest posterior density credible set.  \citet{reimagined} discusses cases in which $\prior_x^\star$ has a density function that satisfies the above comonotonicity requirement.  In the present case where $\prior_x^\star$ has both discrete and continuous components, there is no ``density'' in the usual sense, hence the focus on more general credible set constructions.  One obvious choice of $f$ that satisfies the conditions of Theorem~\ref{thm:coverage} is $f_x = \pi_x$, and the relevant coverage probability properties follow immediately from the defining property of $\prior_x^\star$, namely, \eqref{eq:char}.  Another reasonable strategy is to take $f_x(\theta) = R(x,\theta)$, the relative likelihood.  Generally, this $f_x$ will not be comonotone with $\pi_x$, so the conclusions of Theorem~\ref{thm:coverage} do not apply, at least not exactly.  In our experience, however, this does tend to work well in the following sense.  Define the function 
\[ U(x,\theta) = \prior_x^\star(\{\vartheta: R(x,\vartheta) \leq R(x,\theta)\}, \quad \theta \in \TT. \]
Then the credible set described above, with $f_x(\theta)=R(x,\theta)$, would achieve the nominal frequentist coverage probability, if the random variable $U(X,\Theta)$, as a function of $X \sim \prob_\Theta$, is stochastically no smaller than $\unif(0,1)$.  Figure~\ref{fig:coverage} shows distribution functions of $U(X,\Theta)$ for three different sample sizes, each based on 1000 samples from the binomial model with $\Theta=0.6$.  If the distribution function is on or below the diagonal, $\unif(0,1)$ reference line, then the credible set achieves the frequentist's nominal coverage at that level.  In this case, the distribution function is below the reference line, at least in the lower left-hand corner near the typical levels like 0.05 and 0.10, which explains our claim that credible sets based on the relative likelihood ``tend to work well.'' 

A natural question is if the cumulative distribution function, $Q_x^\star$, derived from the measure $\prior_x^\star$ can be used to construct a credible set with nominal coverage?  For example, the equal-tailed credible interval whose endpoints are obtained by setting $Q_x(\theta)$ equal to $\alpha/2$ and $1-\alpha/2$, respectively, and solving for $\theta$.  In general, such distribution function-based credible intervals are not confidence intervals, but a better understanding the distribution function will be useful shortly.  For some simple insights, start with the case of a smooth contour $\theta \mapsto \pi_x(\theta)$ as considered in \citet{reimagined} and in Appendix~\ref{SS:smooth} below.  In particular, the contour is assumed to be non-decreasing from 0 to 1 on the left of the mode, which is the model's maximum likelihood estimator $\hat\theta_x$, and then non-increasing from 1 to 0 to the right of the mode. In this case, the corresponding distribution function has a simple two-part formula:
\begin{equation}
\label{eq:cdf}
Q_x^\star(\theta) \equiv \prior_x^\star\bigl( (-\infty, \theta] \bigr) = \begin{cases} \frac12 \pi_x(\theta) & \text{if $\theta \leq \hat\theta_x$} \\ 1 - \frac12 \pi_x(\theta) & \text{if $\theta > \hat\theta_x$}. \end{cases} 
\end{equation}
Two key things about the smooth case are clear from this expression.  First, the credible interval obtained by setting the distribution function equal to $\alpha/2$ and $1-\alpha/2$, respectively, and solving equations is exactly the contour's $\alpha$-cut, which is a genuine confidence interval according to \eqref{eq:coverage}.  Therefore, in the smooth case, the distribution function-based credible interval derived from the probabilistic approximation of a possibilistic IM is indeed a confidence interval.  Second, in the smooth case, it is clear from the above display that the mapping $\mu$ that sends the contour function $\pi_x(\cdot)$ to the corresponding distribution function $Q_x(\cdot)$ associated with its inner probabilistic approximation is continuous with respect to the uniform topology.  That is, if $\pi_x$ and $\pi_x'$ are two smooth possibility contours that differ uniformly by no more than $\eps$, then 
\[ \bigl| \mu(\pi_x)(\theta) - \mu(\pi_x')(\theta) \bigr| = \begin{cases} \bigl| \frac12\{ \pi_x(\theta) + \pi_x'(\theta)\} - 1 \bigr| & \text{$\theta$ between the two modes} \\ \frac12|\pi_x(\theta) - \pi_x'(\theta)| & \text{otherwise}, \end{cases} \]
and, clearly, the right-hand side is no more than $\eps/2$, so $\mu$ is continuous.  

Unfortunately, there is no simple formula for $Q_x^\star$ when the contour has jump discontinuities.  In fact, it is very easy to see, e.g., using the case presented in Figure~\ref{fig:jumpy}, that the distribution function is generally different from $\frac12 \pi_x(\theta)$ to the left of its mode.  In an extreme case, a jump of height 1 on the right side of the mode results in $Q_x^\star(\theta) = \pi(\theta)$, which is quite different from the expression in \eqref{eq:cdf}.  All that can be said is that, if $\theta$ is to the left of the mode, then the defining property \eqref{eq:char} of the inner probabilistic approximation $\prior_x^\star$ can be re-expressed in terms of its distribution function $Q_x^\star$:
\[ \pi_x(\theta) = Q_x^\star(\theta) + \{1 - Q_x^\star(\pi_x^{-R}(\theta))\}, \]
where $\pi_x^{-R}(\theta)$ solves the equation $\pi_x(t) = \pi_x(\theta)$ for $t$ to the right of the mode.  An important observation, however, is that the distribution function formula \eqref{eq:cdf} for the case of a smooth contour also holds in the not-smooth case, provided that the contour is symmetric.  Indeed, when the contour is symmetric, there are no shadow intervals, so any change in the contour---whether smooth or jumpy---is reflected equally on both sides of the mode.  So, in the symmetric case, the two terms on the right-hand side of the above display are equal and, consequently, $Q_x^\star(\theta) = \frac12 \pi_x(\theta)$, just like in \eqref{eq:cdf}.  A similar analysis can be given for the case when $\theta$ is to the right of the mode.  

While we cannot expect the contour to be exactly symmetric, we can expect approximate symmetry like what is observed in Figure~\ref{fig:toy}(b). A bit more formally, the possibilistic Bernstein--von Mises theorem from \citet{imbvm.ext} shows that, under the usual regularity conditions, the contour will, for large sample sizes, resemble a Gaussian contour and, therefore, will be approximately symmetric.  Even more formally, define a version of the IM's possibility contour in the local coordinate system, i.e., 
\[ \check \pi_n(z) = \pi_{x^n}(\hat\theta_{x^n} + n^{-1/2} z), \quad z \in \RR. \]
\citet{imbvm.ext} prove that, for any fixed compact $\mathcal{K} \subset \RR$, 
\[ \sup_{z \in \mathcal{K}} |\check\pi_n(z) - \gamma(z)| \to 0, \quad \text{in $\prob_\Theta$-probability as $n \to \infty$}, \]
where $\gamma=\gamma_\Theta$ is the Gaussian contour with mean 0 and variance $I_\Theta^{-1}$, i.e., 
\[ \gamma(z) = 1 - F(I_\Theta \, z^2), \]
with $F$ the $\chisq(1)$ distribution function and $I_\Theta$ the model's Fisher information evaluated at the true $\Theta$. As is clear from the definition, the Gaussian contour $z \mapsto \gamma(z)$ is symmetric around 0.  Of course, those regularity conditions are satisfied for binomial, Poisson, and many other discrete-data models that we have in mind here.  It should come as no surprise that the possibilistic Bernstein--von Mises theorem implies a version of the classical/probabilistic Bernstein--von Mises theorem for the corresponding inner probabilistic approximation.  The following theorem establishes this formally, extending one of the main results in \citet{reimagined} to cases with jump discontinuities.  A relevant technical point is that the ``total-variation convergence'' in classical Bernstein--von Mises theorems is out of reach because the suitably scaled $\prior_x^\star$ and the Gaussian limit do not share the same dominating measure; so, ``weak convergence'' is the natural target. 

\begin{thm}
\label{thm:bvm}
Under the ``usual regularity conditions''---see Appendix~\ref{AA:proofs}---the centered and scaled inner probabilistic approximation $\check\prior_{X^n}^\star$ of the possibilistic IM $\uPi_{X^n}$ converges weakly to $\nm(0, I_\Theta^{-1})$ in $\prob_\Theta$-probability as $n \to \infty$. 
\end{thm}

\begin{figure}[t]
\begin{center}
\subfigure[$n=20$]{\scalebox{0.39}{\includegraphics{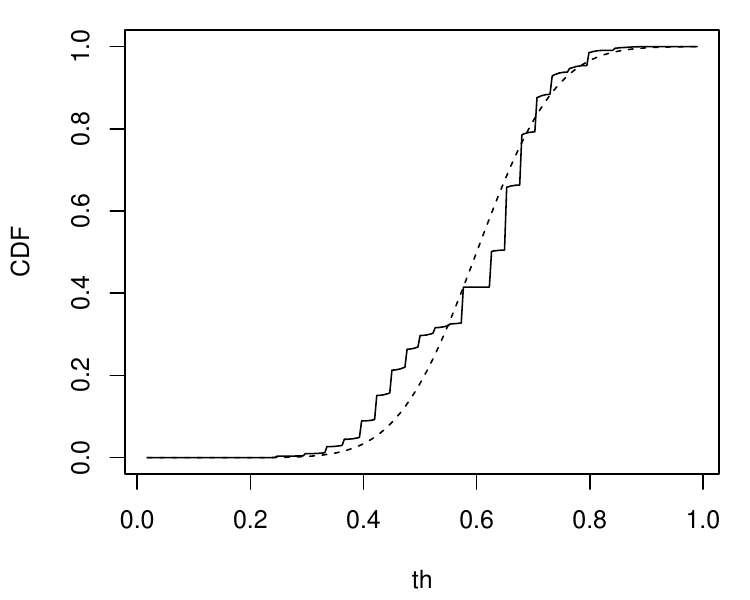}}}
\subfigure[$n=40$]{\scalebox{0.39}{\includegraphics{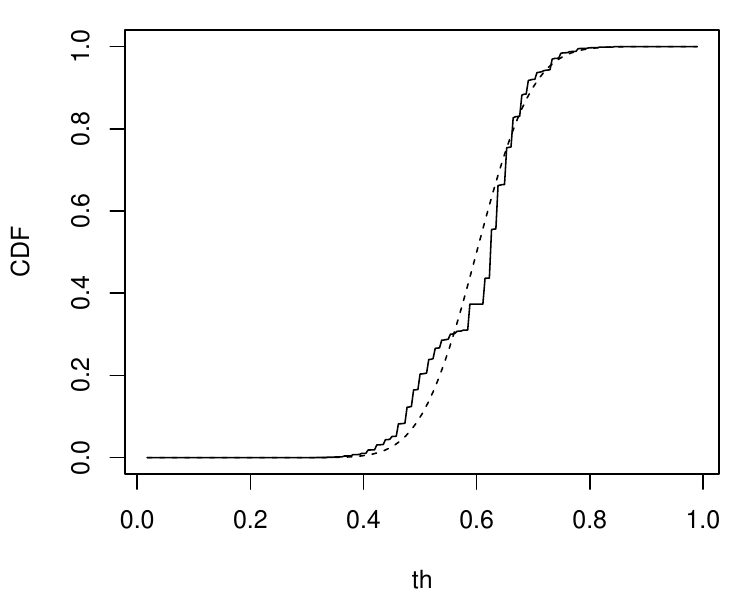}}}
\subfigure[$n=80$]{\scalebox{0.39}{\includegraphics{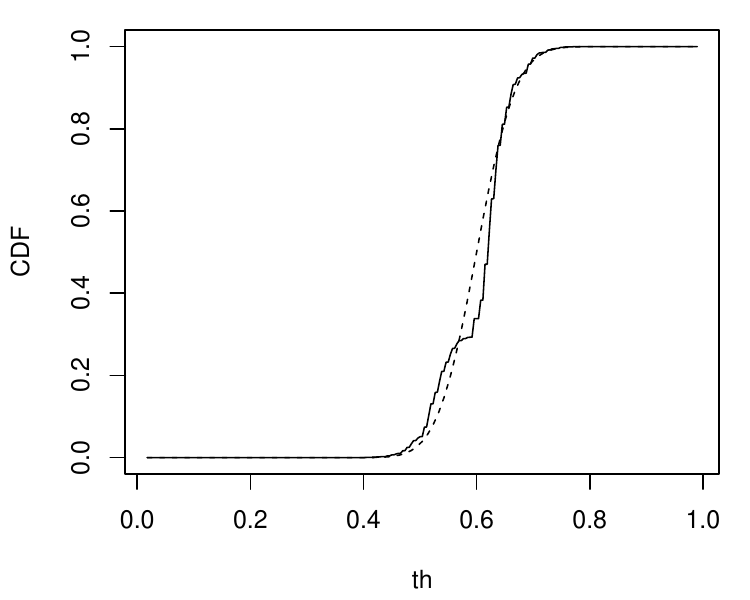}}}
\end{center}
\caption{Plots of the distribution function of the inner probabilistic approximation (solid) compared to the Gaussian limiting approximation (dashed) in the binomial model based on three data sets: all three have $\hat\theta = 0.6$, with varying sample sizes $n$.}
\label{fig:cdf}
\end{figure}

As a quick numerical illustration, Figure~\ref{fig:cdf} shows the distribution functions $Q_x^\star$ associated with the inner probabilistic approximation corresponding to the discrete and continuous parts as displayed in Figure~\ref{fig:discrete} using the same data sets.  Visually, the approximation is already quite accurate even for relatively low sample sizes, which corroborates the result presented in Theorem~\ref{thm:bvm}.  Perhaps unexpectedly, the Gaussian approximation is worst near the maximum likelihood estimator, which is a consequence of the contour's plateau there and $\prior_x^\star$'s choice to push the mass to the tails as much as possible to ensure calibration.  Of course, the discrepancy between the two distribution functions will turn invisible/negligible as the sample size increases further.  

The key take-away message is that the IM's inner probabilistic approximation achieves the same asymptotic optimality expected of fiducial or Bayesian posterior distributions without any specialized argument or choice of default prior distribution.  It is in this sense that our proposal serves as a ``reIMagined'' fiducial argument or no-prior Bayes solution, one that achieves certain finite-sample calibration guarantees (e.g., Theorem~\ref{thm:coverage}) without sacrificing on asymptotic efficiency.

\section{Conclusion}
\label{S:discuss}

Following up on the recent statistical developments in \citet{immc, reimagined} related to probabilistic approximations of possibility measures, where we offer more explicit formulas for those approximations in several different cases: the main paper focuses on the most complicated case where the possibility contour has jump discontinuities, whereas two other related cases are presented in Appendix~\ref{A:inner}.  The key insight in the case with jump discontinuities is that the probabilistic approximation will typically have atoms, i.e., it will be a mixture of discrete and continuous distributions, and a strategy for identifying both the discrete and continuous components.  In the case of possibilistic IMs, these probabilistic approximations can be viewed as alternatives to the more familiar default-prior Bayes and fiducial solutions.  In the challenging case where the data are discrete, the possibilistic IM's contour has jump discontinuities and, hence, the corresponding probabilistic approximation has atoms---this is unique to our new approach.  We show that the corresponding probability distribution has a number of desirable properties, including an exact calibration property which facilitates the construction of credible sets with exact frequentist coverage as well as a version of the Bernstein--von Mises theorem.  

While many of the insights presented here apply for general discrete-data models, the details presented here focus on the scalar parameter case.  What makes the extension to multi-parameters cases particularly challenging is that the jump discontinuities that were critical to the characterization in the scalar case above are no longer just isolated points: in the two-dimensional case, for example, the contour resembles a mountain with jagged cliffs, so the jump discontinuities can be described as a union of one-dimensional subsets in the two-dimensional parameter space.  Work is underway to extend the results in the present paper to the general multi-parameter discrete data models like multinomial.  In a forthcoming follow-up work, we have developed alternative strategies for carrying out possibilistic IM computations using probabilistic approximations, which is able to accurately handle categorical regression models, among others.

\section*{Acknowledgments}

This work is partially supported by the U.S.~National Science Foundation, DMS--2412628, and the National Heart, Lung, And Blood Institute of the National Institutes of Health, T32HL079896.

\appendix

\section{Technical details from the main paper}
\label{A:details}

\subsection{More background on possibility theory}
\label{AA:background}

The main paper is largely focused on finding a probabilistic approximation of a given possibility measure.  But we think it might be helpful to supplement that emphasis with a few details about the opposite case: where we seek a possibilistic approximation of a given probability measure.  The point is that the former approximation is, in a sense to be made clear below, the inverse of the latter approximation.  

Consider a generic case with an uncertain variable $T \in \TT$.  Suppose that uncertainty about $T$ is quantified using a probability distribution $\prob$.  What is the ``best approximation'' of the probability distribution $\prob$ by a possibility $\uPi$?  Since the possibility measure $\uPi$ has a credal et $\cred(\uPi)$, it makes sense to choose $\uPi$ in such a way that $\cred(\uPi)$ is the smallest of all those credal sets that contains the probability $\prob$.  With few exceptions, this is not enough information to pin down on a specific possibilistic approximation.  But given a suitable ranking function---or {\em elementary plausibility ordering} \citep{hose2022thesis}---a unique solution can be found.  Indeed, given $\rho$, the possibility $\uPi=\uPi^\rho$ with contour 
\begin{equation}
\label{eq:prob.to.poss}
\pi(t) = \pi^\rho(t) = \prob\{ \rho(T) \leq \rho(t) \}, \quad t \in \TT, 
\end{equation}
has the following property: among all those possibility measures whose contours are comonotone with $\rho$ and whose credal sets contain $\prob$, the possibility $\uPi=\uPi^\rho$ is such that $\cred(\uPi^\rho)$ is a subset of all these other credal sets.  When $\prob$ has a density or mass function $p$, it is common to take $\rho=p$; in that case, the mapping from $\prob$ to $\uPi$ defined in \eqref{eq:prob.to.poss} is often called a {\em probability-to-possibility transform}.  

Looking at \eqref{eq:char} from this perspective, it is not difficult to see that, in the main paper, the distribution $\prior^\star$ that we seek is one such that the probability-to-possibility transform of $\prior^\star$ with respect to ranking $\rho=\pi$ is exactly equal to the possibility measure with contour $\pi$.  In other words, $\prior^\star$ is the {\em possibility-to-probability} transform, or the inverse probability-to-possibility transform of $\uPi$, again, relative to the ranking $\rho=\pi$.

\subsection{More background on possibilistic IMs}
\label{AA:background2}

A detail relevant to the developments started in Section~\ref{SS:approx} is the claim that replacing the contour $\pi=\pi_x$ by its least monotone majorant is ``inconsequential.''  Next we provide a detailed justification of this claim.  

As shown in \citet{imbvm.ext}, the contour will converge uniformly (on compacts) to the Gaussian contour, which is monotone in the sense we're considering here.  Therefore, the problematic non-monotonicity that appears in examples can't be extreme in the sense that, if the sample size is sufficiently large, then the difference between the contour and its least monotone majorant is very small.  To make this more formal, consider comparing the sup-norm distance between a contour $\pi$ and the Gaussian contour $\gamma$ and between the monotone majorant $\pi^\text{\sc maj}$ and $\gamma$.  Without loss of generality, imagine an interval between jumps where $\gamma$ is strictly increasing, $\pi$ is non-monotone, and $\pi^\text{\sc maj}$ is farther from $\gamma$ on that interval than $\pi$; in that case, $\pi$ must be over-estimating $\gamma$. But the distance from $\gamma$ is largest for both $\pi$ and $\pi^\text{\sc maj}$ at the jump on left endpoint of this interval, and those two distances are equal.  Therefore, the supremum distance can't be any greater for $\pi^\text{\sc maj}$ than it is for $\pi$, which means that $\pi^\text{\sc maj}$ converges uniformly to $\gamma$.  So, at least when the sample size is moderate to large, anything that can be said about $\pi$ also applies to $\pi^\text{\sc maj}$, and it's in this sense that we say that the change from $\pi$ to its least monotone majorant is ``inconsequential.''

\subsection{Proofs}
\label{AA:proofs}

\subsubsection{Proof of Theorem~\ref{thm:jumpy}}

\begin{figure}[t]
\begin{center}
\scalebox{0.55}{\includegraphics{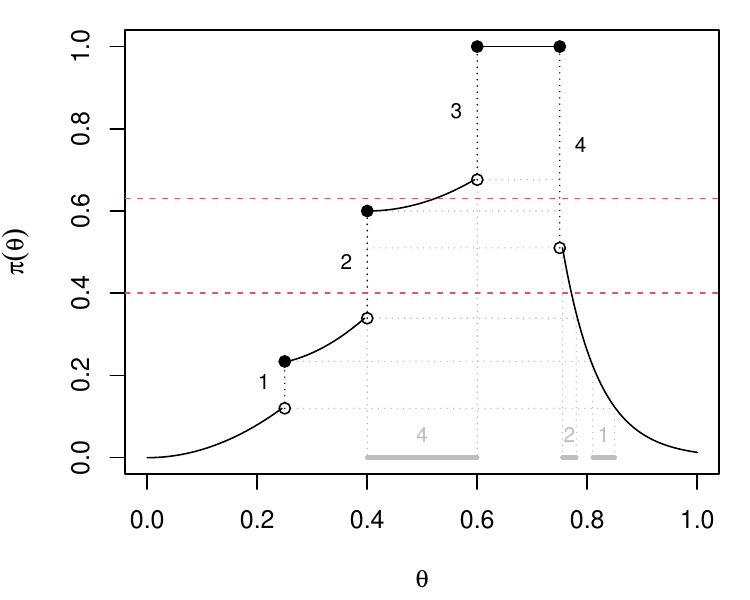}}
\end{center}
\caption{Slightly modified version of Figure~\ref{fig:jumpy} in the main text.  The only difference is the inclusion of the horizontal red dashed lines referred to in the proof of Theorem~\ref{thm:jumpy}.}
\label{fig:jumpy.proof}
\end{figure}

The goal is to verify \eqref{eq:char}.  Let $\alpha$ be such that there exists $\theta$ with $\pi(\theta)=\alpha$.  The dashed lines (red) in Figure~\ref{fig:jumpy.proof} represent two such $\alpha$ values.  Start with the smaller of those two $\alpha$'s.  Let $(-\infty, \theta_L(\alpha))$ and $(\theta_R(\alpha), \infty)$ denote the intervals where $\pi(\theta) \leq \alpha$ on the left and right, respectively.  Then we need to evaluate the $\prior^\star$-probability of the union of these two disjoint intervals.  Following the case in Figure~\ref{fig:jumpy.proof} as a concrete example, the left-hand side involves integrating the absolutely continuous part over $(-\infty,t_1)$ and $(t_1, t_2)$; but there is nothing to add from the discrete part because the jump at $t_1$ is assigned 0 mass.  So, the $\prior^\star$-probability on the left-hand side of the mode is 
\[ \int_{-\infty}^{t_1} \tfrac12 \dot\pi(\theta) \, d\theta + \int_{t_1}^{t_2} \tfrac12 \dot\pi(\theta) \, d\theta = \tfrac12 \{ \pi(t_1-) + \pi(t_2-) - \pi(t_1+) \}. \]
Since this range contains none of shadow intervals, the factor $\frac12$ appears throughout.  Next, on the right-hand side, there are no jumps in the relevant range so we only need to integrate the absolutely continuous part from $\theta_R(\alpha)$ to infinity.  But there are two shadow intervals, $\mathcal{T}_1$ and $\mathcal{T}_2$, that we have to contend with, where $\theta_R(\alpha)$ is contained in $\mathcal{T}_2$.  Ultimately, that right-hand side integral boils down to
\[ \int_{\theta_R(\alpha)}^{\max \mathcal{T}_2} |\dot\pi(\theta)| \, d\theta + \int_{\max \mathcal{T}_2}^{\min \mathcal{T}_1} \tfrac12 |\dot\pi(\theta)| \, d\theta + \int_{\mathcal{T}_1} |\dot\pi(\theta)| \, d\theta + \int_{\max \mathcal{T}_1}^\infty \tfrac12 |\dot\pi(\theta)| \, d\theta, \]
which equals
\begin{align*}
\pi(\theta_R(\alpha)) & - \pi(\max\mathcal{T}_2) + \tfrac12\{ \pi(\max\mathcal{T}_2) - \pi(\min\mathcal{T}_1)\} \\
& \qquad + \pi(\min\mathcal{T}_1) - \pi(\max\mathcal{T}_1) + \tfrac12 \pi(\max\mathcal{T}_1) \\
& = \pi(\theta_R(\alpha)) - \tfrac12 \{ \pi(\max\mathcal{T}_2) + \pi(\max\mathcal{T}_1) - \pi(\min\mathcal{T}_1) \} \\
& = \alpha - \tfrac12 \{ \pi(t_2-) + \pi(t_1-) - \pi(t_1+)\}.
\end{align*}
Adding this to the left's $\prior^\star$-probability derived above, we get the desired result:
\[ \prior^\star(\{\theta: \pi(\theta) \leq \alpha\}) = \alpha. \]
To further illustrate the nuance involved in this construction and derivation, we repeat the above analysis but for the other larger $\alpha$ value shown in Figure~\ref{fig:jumpy.proof}.  Starting with the right-hand side in this case, we integrate the absolutely continuous part from $t_4-$ to $\infty$.  We can easily do this by just continuing the latter calculation in the above display, replacing the $\theta_R(\alpha)$ there with $t_4-$, so the $\prior^\star$-probability assigned to the right is 
\[ \pi(t_4-) - \tfrac12 \{ \pi(t_2-) + \pi(t_1-) - \pi(t_1+)\}. \]
On the left-hand side, in addition to the part 
\[ \tfrac12 \{ \pi(t_1-) + \pi(t_2-) - \pi(t_1+) \} \]
from the previous calculation, we add the discrete portion 
\[ \pi(t_2) - \pi(t_4-) \]
associated with the jump at $t_2$ and the continuous portion 
\[ \int_{t_2+}^{\theta_L(\alpha)} q_\text{\sc c}^\star(\theta) \, d\theta = \int_{t_2+}^{\theta_L(\alpha)} \dot\pi(\theta) \, d\theta = \alpha - \pi(t_2). \]
Note that the factor $\frac12$ is missing in the above calculation because this range of integration belongs to a shadow interval $\mathcal{T}_4$.  Putting each of these pieces together gives the desired result, $\prior^\star(\{\theta: \pi(\theta) \leq \alpha\}) = \alpha$, as in \eqref{eq:char}.

\subsubsection{Proof of Theorem~\ref{thm:discrete.ipa}}

The claim follows from Theorem~\ref{thm:jumpy}, once we demonstrate that the possibilistic IM's contour function has a pair of jumps creating a hole that goes through both sides of the contour, like how the jumps labeled 2 and 3 in Figure~\ref{fig:jumpy} overlap with the jump labeled 4.  First we establish the existence of jumps, which we think is an insightful detail, then show that there are overlapping jumps in the sense just described.  Towards this, for a given data set $x$, we prove that the IM contour $\pi_x$ has a jump at $\theta$ if there exists another data set $y$ such that $R(y,\theta) = R(x,\theta)$.  Recall that the possibility contour is a sum:
\begin{equation}
\label{eq:pi.sum}
\pi_x(\theta) = \sum_{z \in \mathcal{R}(x,\theta)} p_\theta(z), 
\end{equation}
where $\mathcal{R}(x,\theta) = \{z: R(z,\theta) \leq R(x,\theta)\}$.  Without loss of generality, suppose that $R(y,\theta+\eps) > R(x,\theta)$, for some small $\eps > 0$.  Then $y \not\in \mathcal{R}(x,\theta + \eps)$ and so $y$ does not contribute to the sum that determines $\pi_x(\theta + \eps)$.  But once $\eps = 0$, and not before, the data set $y$ is included in $\mathcal{R}(x,\theta)$ and, hence, contributes an amount $p_\theta(y)$ to the sum that determines $\pi_x(\theta)$.  It is this contribution that creates the jump.  It remains to show that there exists a $y$ such that $R(y,\theta) = R(x,\theta)$, but this is actually trivial.  To simplify the discussion, without loss of generality, identify $x$ and generic data set $y$ with the corresponding values of their minimal sufficient statistic; the relative likelihood only depends on the minimal sufficient statistic, so this value is all that really matters.  Then every minimal sufficient statistic $y$ that differs from $x$ admits such a solution $\theta$.  Indeed, for each such $y$, the function $\theta \mapsto R(y,\theta) - R(x,\theta)$ takes a positive value at $\theta=\hat\theta_y$ and negative value at $\theta=\hat\theta_x$.  Since this function is continuous, the intermediate value theorem implies that it has a root, at which we showed that $\pi_x$ has a jump.   

Next, we establish existence of a {\em special} pair of overlapping jumps that creates a hole as described above.  This suffices show that the probabilistic approximation has an atom and, hence, prove the claim, but this does not describe the discrete component of $\prior_x^\star$ in its entirety---in general, there will be multiple atoms, as can be readily seen just by looking at the plot of the contour, like in Figure~\ref{fig:toy}(b).  Recall that, by construction, $\pi_x(\hat\theta_x)=1$.  Let $y^\text{\sc r}$ and $y^\text{\sc l}$ denote the two data sets, distinct from $x$, such that $\hat\theta_{y^\text{\sc r}}$ and $\hat\theta_{y^\text{\sc l}}$ are as close to, but different from, $\hat\theta_x$ on the right and left, respectively.  That is, 
\[ \hat\theta_{y^\text{\sc l}} < \hat\theta_x < \hat\theta_{y^\text{\sc r}}, \]
and that there is no other $y$ with $\hat\theta_y$ closer to $\hat\theta_x$ than these two bounds.  If $\hat\theta_x$ is on a boundary of the parameter space, so that one can only find $y^\text{\sc r}$ or $y^\text{\sc l}$ but not both, then the following argument can be easily modified, as we explain below.  The result proved in the above paragraph implies that there are jumps in the contour at $\theta^\text{\sc l}$ and $\hat\theta^\text{\sc r}$ between $\hat\theta_x$ and the two endpoints in the above display, respectively.  As $\theta$ moves away from $\hat\theta_x$ in either direction, the contour drops from 1 to some value strictly less than 1 at each of $\theta^\text{\sc l}$ and $\theta^\text{\sc r}$.  Those two drops, which both originate at value 1, obviously are overlapping, hence there is a hole at the first drop in either direction of the maximum likelihood estimator, which implies existence of atoms in the probabilistic approximation.  In the case mentioned above, where $\hat\theta_x$ is on the boundary of $\TT$, a jump on one side of the mode is enough to create a hole and hence an atom.

\subsubsection{Proof of Theorem~\ref{thm:coverage}}

If $f_x$ is comonotone with $\pi_x$, i..e, if $f_x = g(\pi_x)$ for some increasing function $g$, then both $f_x$ and $\pi_x$ have the same level sets.  More formally, 
\[ f_x(\theta) > c \iff \pi_x(\theta) > g^{-1}(c), \]
so, by \eqref{eq:char}, the cutoff $c_\alpha^f(x)$ defined in Theorem~\ref{thm:coverage} must satisfy 
\[ g^{-1}\{ c_\alpha^f(x) \} = \alpha. \]
Then the credible set $C_\alpha^f(x)$ equals $\{\theta: \pi_x(\theta) > \alpha\}$ and the coverage probability claim follows from the key result \eqref{eq:coverage} of the IM construction.

\subsubsection{Proof of Theorem~\ref{thm:bvm}}

A precise statement of the regularity conditions imposed for the possibilistic Bernstein--von Mises theorem established by \citet{imbvm.ext} can be found in that paper.  Here we give just a brief summary of those.  First, it is assumed that the model's parameter space $\TT$ is open and that each point in the interior is ``regular'' in the sense that the standard local asymptotic normality results \citep[e.g.,][Ch.~7]{vaart1998} hold there; roughly, this implies differentiability of the log-likelihood and suggests that a Taylor approximation can be used as in the standard textbook treatments \citep[e.g.,][Theorem~3.10]{lehmann.casella.1998} on maximum likelihood estimators' asymptotic normality.  In addition, there must exist a non-singular Fisher information matrix and certain Lipschitz conditions on the likelihood function that allows for proper handling of the remainder term in the aforementioned Taylor approximation.  Finally, \citet{imbvm.ext} work under the assumption that the maximum likelihood estimator is consistent, which holds in many standard cases, but can be restrictive in other.  We expect that this assumption can be relaxed, but we will not pursue this matter here.  

Now to the proof of Theorem~\ref{thm:bvm}. To prove weak convergence, it is sufficient to show that the distribution function associated with $\check\prior_{X^n}^\star$ converges pointwise, in $\prob_\Theta$-probability, to the corresponding $\nm(0, I_\Theta^{-1})$ distribution function.  As explained in the main text prior to the theorem statement, when the contour is symmetric, the mapping $\mu$ from said contour to the distribution function of its corresponding inner probabilistic approximation is incredibly simple and clearly continuous with respect to the topology induced by uniform convergence on compacts.  Since $\check\pi_n$ converges in $\prob_\Theta$-probability to $\gamma$ in the aforementioned topology, and the symmetric limit $\gamma$ belongs to the continuity set of $\mu$, the desired claim follows from the continuous mapping theorem.  That is, for $z \leq 0$, 
\[ \check Q_x^\star(z) = \mu(\check \pi_x)(z) \to \mu(\gamma)(z) = \Phi(I_\Theta \, z), \quad \text{in $\prob_\Theta$-probability}, \]
with a similar statement holding for $z > 0$.  Since this holds for all $z$, we conclude that $\prior_x^\star$ converges weakly to $\nm(0, I_\Theta^{-1})$ in $\prob_\Theta$-probability, proving the claim.


\section{Probabilistic approximation results}
\label{A:inner}

\subsection{Common thread}
\label{SS:common}

While the specific form of the inner probabilistic approximations vary in the different cases considered below, there is a common thread that links them and informs the general construction.  Recall that the to-be-approximated possibility measure has a contour $\pi$, that maps a space $\TT$ to the interval $[0,1]$; while the contour may not take all the values in $[0,1]$, it does effectively reach 1.  Related to the condition in \eqref{eq:char}, consider the lower $\alpha$-level sets $\{\theta \in \TT: \pi(\theta) \leq \alpha\}$ indexed by $\alpha \in [0,1]$.  When $\alpha=0$, the set is empty, except for those ``impossible'' $\theta$ values with $\pi(\theta)=0$.  Of course, as $\alpha$ increases, more and more points enter into the lower level sets.  How the inner probabilistic approximation assigns mass to $\TT$ is determined by how much $\alpha$ must increase so that a point is just excluded from the lower level set to just included.  

In the case of finite $\TT$, the probability mass assigned to a particular $\theta \in \TT$ is determined by how much $\alpha$ increases between the previous entry into the lower level set and when $\theta$ first enters---this is just the difference between the possibility assigned to $\theta$ and the next smallest possibility value.  If there are ties in the sense that two or more points have the same possibility value, then they enter the lower level set simultaneously, and the corresponding probability mass is shared among all of them.  

Similarly, in the case where $\TT$ is a continuum and $\pi$ is continuous, the probability is assigned in an analogous way: an infinitesimal increase in $\alpha$ introduces new $\theta$ values into the lower level set and the density assigned at a given $\theta$ value is proportional to the absolute rate of change in the contour $\pi$ there.  One subtlety in the continuous case is that there are always ties in the sense above, so the density is shared among all of them. 

The case where the space $\TT$ is a continuum but the contour has jump discontinuities is more complicated to describe, but still corresponds to considering the change in $\alpha$ required for a $\theta$ that is absent from the lower level set to enter.  The subtlety lies largely in the fact that there will often be ties where the contour is smooth, but not where the contour has jumps, so deciding on when to share the mass/density requires some care.

\subsection{Case: finite/countable spaces}
\label{SS:finite}

Consider a finite space $\TT = \{t_1, t_2, \ldots, t_S\}$ and a possibility contour $\pi: \TT \to [0,1]$ defined there.  For simplicity, assume that $\pi$ assigns unique values to each member of $\TT$, i.e., there are no ties.  Generalizations to this relatively simple case will be addressed later.  

First, define the ordered values $t_{(s)}$, for $s=1,\ldots,S$, such that 
\[ \pi(t_{(1)}) < \pi(t_{(2)}) < \cdots < \pi(t_{(S)}) \equiv 1. \]
Next, define a probability measure $\prior^\star$ on $\TT$ with mass function $q^\star$ according to the rule
\begin{align*}
q^\star(t_{(1)}) & = \pi_{(1)} \\
q^\star(t_{(s)}) & = \pi(t_{(s)}) - \pi(t_{(s-1)}),  \quad s= 2,\ldots,S.
\end{align*}
A consequence of the result shown in \citet[][Sec.~2.3]{dubois.etal.2004} is that $\prior^\star$ is a probabilistic approximation of $\uPi$; in fact, when there are no ties, it is the unique such probability distribution.  This is formally stated and proved here for completeness. 

\begin{thm}
\label{thm:finite}
Let $\uPi$ be a possibility measure defined on the finite space $\TT$.  Then the probability measure $\prior^\star$ defined above is an inner probabilistic approximation of $\uPi$.
\end{thm}

\begin{proof}
Let $T$ be a $\TT$-valued random variable with distribution $\prior^\star$.  For any $r \in \{1,\ldots,S\}$, 
\begin{align*}
\prior^\star\{ \pi(T) \leq \pi(t_{(r)}) \} & = \sum_{s=1}^S 1\{ \pi(t_{(s)}) \leq \pi(t_{(r)}) \} \, q^\star(t_{(s)}) \\
& = \sum_{s=1}^r q^\star(t_{(s)}) \\
& = \pi(t_{(1)}) + \sum_{s=2}^r \{ \pi(t_{(s)}) - \pi(t_{(s-1)}) \} \\
& = \pi(t_{(r)}).
\end{align*}
Since any $t \in \TT$ corresponds to $t_{(r)}$ for some $r$, we get that $\prior^\star\{ \pi(T) \leq \pi(t) \} = \pi(t)$ for each $t \in \TT$.  Therefore, $\prior^\star$ is the inner probabilistic approximation of $\uPi$, as claimed. 
\end{proof}

A couple of relevant generalizations are possible.  First, if there are ties in the possibility assignments, then one can effectively treat the set of all those $t$ that share a common $\pi$ value as a single ``point'' and apply the above construction.  Then it's just a matter of splitting the mass assigned to such composite ``points'' among the genuine points that it contains, and this can be done in any manner.  A natural choice is to split the mass equally among the ties so, in case there are two points with the same possibility contour value, the difference defined above would be cut in half.  

Second, if $\TT$ is infinite but countable, then the difficulty is in sorting the $t$ values so that the corresponding $\pi$ values are decreasing---in the infinite case, we may not be able to identify the ``next smallest'' value.  But if the ordering of the $t$ values by $\pi$ is given, e.g., if we know that $\pi$ is decreasing along a known ordering of the $t$'s, then we can define the probabilistic approximation via differences in virtually the same way as above.  


\subsection{Case: smooth contours}
\label{SS:smooth}

\citet[][Sec.~3.2]{dubois.etal.2004} establish a formula for the density of the inner probabilistic approximation of a possibility measure $\uPi$ with contour $\pi$ in terms of the derivative of $\pi$.  Their analysis only covers the case where $\TT$ is a subset of the real line.  Here we extend their result to the $d$-dimensional case, with $d \geq 1$, where $\TT \subseteq \RR^d$.  This result is probably already known, but since it is relevant and we are not aware of any reference where it has been established, we flesh out the details here. 

Let $\pi: \TT \to [0,1]$ be a possibility contour with a unique mode $\tau$, so that $\pi(\tau) = 1$.  Assume that $\pi$ is differentiable everywhere except perhaps at $\tau$, where the peak is often sharp; see Figure~\ref{fig:toy}(a).  Furthermore, assume that $\pi$ is strictly decreasing moving away from the mode, with $\pi(\theta) \to 0$ as $\theta \to \infty$ in any direction.  The proof makes use of the spherical coordinate system, where $\tau$ is the center or anchor, $\rho$ is a radius, and $\delta$ is a direction---a unit vector with $\|\delta\| = ( \delta^\top \delta )^{1/2}=1$.  That is, a generic point $\theta \in \TT$ can be uniquely expressed as $\theta = \tau + \rho \, \delta$, where, of course, $\rho = \|\theta - \tau\|$ and $\delta = (\theta-\tau) / \|\theta - \tau\|$.  Note that $\theta$ is $d$-dimensional and so is the corresponding $(\rho, \delta)$, since the unit vector constraint makes $\delta$ $(d-1)$-dimensional.  Define the function
\[ \varpi_\delta(\rho) = \pi(\tau + \rho \, \delta), \quad \rho > 0, \]
indexed by a direction $\delta$.  Then $\rho \mapsto \varpi_\delta(\rho)$ is strictly decreasing in $\rho$ for a given $\delta$, and hence is invertible; it is also differentiable.  

\begin{thm}
\label{thm:smooth}
Under the conditions on $\pi$ described above, define the density 
\begin{equation}
\label{eq:smooth}
q^\star(\theta) = \frac{1}{\sigma(\mathbb{S})} \, \frac{(\theta - \tau)^\top}{\|\theta-\tau\|^d} \, \dot\pi(\theta), \quad \theta \in \TT, 
\end{equation}
where $\mathbb{S}$ is the unit sphere in $\RR^d$, $\sigma(\mathbb{S}) = 2\pi^{d/2}/\Gamma(d/2)$ is the surface area of $\mathbb{S}$, and $\dot\pi$ is the gradient $\pi$.  Then the probability measure $\prior^\star$ determined by the density $q^\star$ above is an inner probabilistic approximation of the possibility measure $\uPi$. 
\end{thm}

\begin{proof}
Let $\prior^\star$ be an inner probabilistic approximation of $\uPi$ and let $\Theta \sim \prior^\star$ be a $\TT$-valued random variable with distribution $\prior^\star$.  Then 
\[ \prior^\star\{ \pi(\Theta) \leq \pi(\theta) \} = \pi(\theta) \quad \text{for all $\theta \in \TT$}. \]
If $q^\star$ is the density corresponding to $\prior^\star$, then the equation in the above display can be rewritten as an integral:
\[ \int_{\vartheta: \pi(\vartheta) \leq \pi(\theta)} q^\star(\vartheta) \, d\vartheta = \pi(\theta). \]
Changing to spherical coordinates as indicated above, the above integral equals 
\[ \int_\mathbb{S} \int_{\varpi_\delta^{-1}(\pi(\theta))}^\infty q^\star(\tau + \rho \, \delta) \, \rho^{d-1} \, d\rho \, J(\delta) \, d\delta, \]
where $\rho^{d-1} \, J(\delta)$ is a Jacobian term associated with the transformation $\vartheta \to \tau + \rho \, \delta$; the form of $J(\delta)$ will not be needed in what follows.  If $q^\star$ is as in \eqref{eq:smooth}, then 
\[ q^\star(\tau + \rho \, \delta) = \frac{1}{\sigma(\mathbb{S})} \, \frac{1}{\rho^{d-1}} \, \delta^\top \dot\pi(\tau + \rho \, \delta) = \frac{1}{\sigma(\mathbb{S})} \, \frac{1}{\rho^{d-1}} \, \dot\varpi_\delta(\rho), \]
and the integral reduces to 
\[ \frac{1}{\sigma(\mathbb{S})} \int_\mathbb{S} \int_{\varpi_\delta^{-1}(\pi(\theta))}^\infty \frac{1}{\rho^{d-1}} \, \dot\varpi_\delta(\rho) \, \rho^{d-1} \, d\rho \, J(\delta) \, d\delta = \frac{1}{\sigma(\mathbb{S})} \int_\mathbb{S} \int_{\varpi_\delta^{-1}(\pi(\theta))}^\infty \dot\varpi_\delta(\rho) \, d\rho \, J(\delta) \, d\delta. \]
By the fundamental theorem of calculus, the inner integral equals $\pi(\theta)$.  Since the outer integral is normalized, the entire display above equals $\pi(\theta)$, as was to be proved. 
\end{proof}

Then the result in \citet[][Sec.~3.2]{dubois.etal.2004} for scalar variables immediately follows as special case of Theorem~\ref{thm:smooth}.

\begin{cor}
For $d=1$, the density of the inner probabilistic approximation is 
\[ q^\star(\theta) = \tfrac12 |\dot\pi(\theta)| = \begin{cases} \tfrac12 \dot\pi(\theta) & \text{if $\theta < \tau$} \\ -\tfrac12 \dot\pi(\theta) & \text{if $\theta > \tau$}, \end{cases} \]
which agrees with the result in \citet[][Sec.~3.2]{dubois.etal.2004}; see, also, \citet{reimagined}.
\end{cor}

There is nothing special about the above case where the mode is in the interior of the space $\TT$.  One can obviously imagine, say, a scalar variable setting in which the contour is strictly increasing or decreasing, so that the mode is on the boundary of the space.  In such cases, it is easy to modify the result in Theorem~\ref{thm:smooth} to accommodate this.  Generalizations to the $d$-dimensional space can be developed, but we will not pursue this here. 

\begin{cor}
\label{cor:monotone}
For $d=1$, if the smooth contour's mode is on the boundary of the space, so that the contour is strictly increasing/decreasing, then the density of the corresponding inner probabilistic approximation is $q^\star(\theta) = |\dot\pi(\theta)|$. 
\end{cor}

Situations where the mode is on the boundary are particular instance of a more general case where the contour does not vanish at the boundary.  Again, focusing on scalar parameters only, this non-vanishing-at-a-boundary property is handled in the same way we handled jump discontinuities in the main text.  That is, the ``jump'' at the boundary casts a shadow on the opposite side of the mode, and that impacts the scale of the density in the shadow interval.  The following corollary considers the case where the contour is smooth except for the jump at the left-most boundary, which casts a shadow $\T$ to the right of the mode. 

\begin{cor}
\label{cor:boundary}
For $d=1$, suppose the smooth contour $\pi$ does not vanish at the left endpoint of its support, which casts a shadow $\T$ interval to the right of the mode.  Then the density of the corresponding inner probabilistic approximation is 
\[ q_x^\star(\theta) = \begin{cases} \frac12 |\dot\pi(\theta)| & \text{if $\theta \not\in \T$} \\ |\dot\pi(\theta)| & \text{if $\theta \in \T$}. \end{cases} \]
\end{cor}


\subsection{Statistical applications}

\begin{figure}[t]
\begin{center}
\scalebox{0.55}{\includegraphics{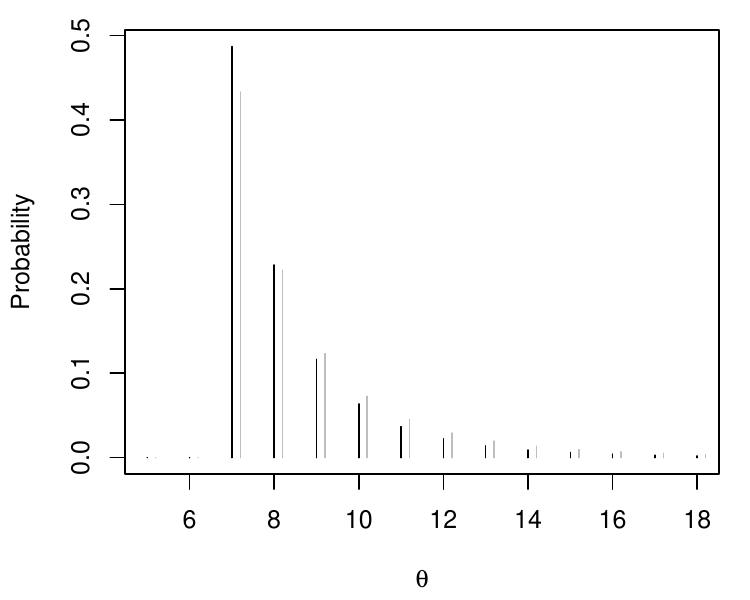}}
\end{center}
\caption{Plot of the mass function corresponding to the probabilistic approximation $\prior_x^\star$ of the possibilistic IM in the discrete uniform example, with $n=5$ and $x_{(n)}=7$. The slightly-offset gray mass function is that of the Bayesian posterior with a flat prior.} 
\label{fig:uniform}
\end{figure}

\begin{ex}
\label{ex:unif.discrete}
Consider a model $\prob_\theta = \unif(\{1,2,\ldots,\theta\})$, where $\theta$ belongs to the countable parameter space $\TT = \{1,2,\ldots\}$.  Then the data $X$ consists of iid samples $X_1,\ldots,X_n$ from $\prob_\Theta$ where $\Theta \in \TT$ is unknown.  While this may be a ``toy problem,'' there are real survey sampling applications \citep[e.g.,][]{basu1969, basu1971}.  Then the relative likelihood is
\[ R(x,\theta) = \Bigl( \frac{x_{(n)}}{\theta} \Bigr)^n \, 1(x_{(n)} \leq \theta), \]
where $x_{(n)}$ is maximum of the observations in the $n$-vector $x$.  This is a strictly decreasing function for $\theta$ on the range $x_{(n)}, x_{(n)}+1, \ldots, N$.  From here, as in \citet{martin.basu}, it is easy to find a closed-form expression for the possibilistic IM's contour function:
\[ \pi_x(\theta) = \prob_\theta\{ R(X,\theta) \leq R(x,\theta) \} = R(x,\theta). \]
This means that the IM contour $\pi_x(\theta)$ is likewise strictly decreasing for $\theta$ bigger than $x_{(n)}$, so we can apply the differencing strategy described in Section~\ref{SS:finite} to get an inner probabilistic approximation $\prior_x^\star$ of $\uPi_x$, which has mass function 
\begin{align*}
q_x^\star(\theta) & = R(x,\theta) - R(x,\theta + 1) \\
& \propto p_\theta(x) - p_{\theta+1}(x) \\
& \propto p_\theta(x) \Bigl\{ 1 - \frac{p_{\theta+1}(x)}{p_\theta(x)} \Bigr\}, \qquad \theta=x_{(n)}, x_{(n)}+1, \ldots. 
\end{align*}
The last two lines in the above display are not necessary to evaluate the mass function, but they do assist in comparing $\prior_x^\star$ to a Bayesian posterior distribution.  Indeed, for $\prior_x^\star$ to be a Bayesian posterior distribution, it is necessary that the last term in the above display, the one in curly braces, be constant in $x$.  Of course, that would require $p_\theta(x)$ to factor as a function of $x$ times a function of $\theta$, which can never happen: the famous Fisher--Neyman factorization theorem says the likelihood always involves a factor that depends directly on both data and parameter.  Therefore, the probabilistic approximation {\em cannot} be a Bayesian posterior distribution for any prior.  A quick numerical illustration---see Figure~\ref{fig:uniform}---with $n=5$ and $x_{(n)}=7$ confirms that $\prior_x^\star$ differs from the Bayesian posterior with a constant/flat prior on $\TT$.
\end{ex}

\begin{figure}[t]
\begin{center}
\subfigure[Possibility contour]{\scalebox{0.39}{\includegraphics{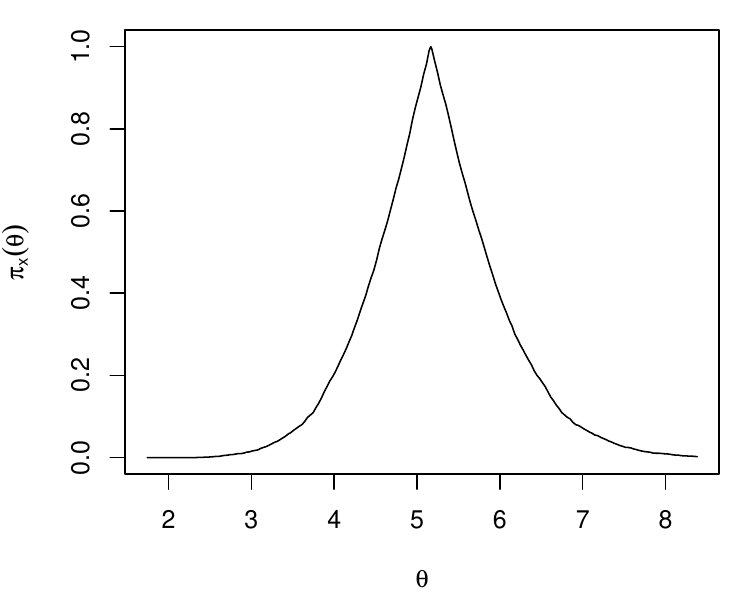}}}
\subfigure[Distribution function]{\scalebox{0.39}{\includegraphics{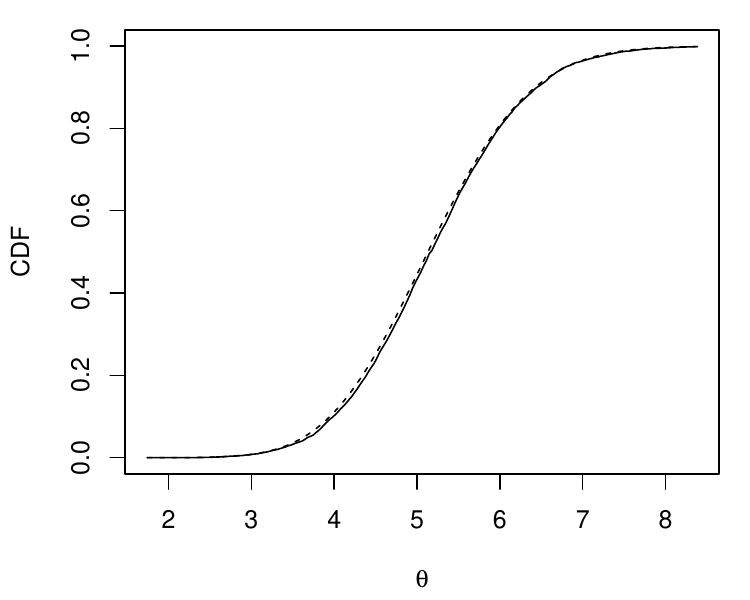}}}
\subfigure[Density function]{\scalebox{0.39}{\includegraphics{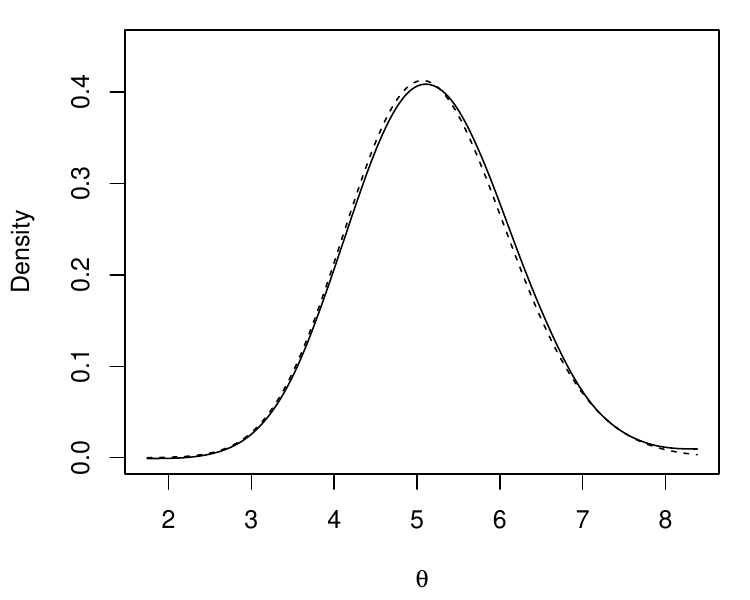}}}
\end{center}
\caption{Results for the gamma illustration in Example~\ref{ex:gam}.  Panels (b) and (c) compare the IM results (solid) with the Jeffreys-prior Bayes results (dashed).}
\label{fig:gam}
\end{figure}

\begin{ex}
\label{ex:gam}
Consider the model $\prob_\theta = \gam(\theta, 1)$, i.e., a gamma distribution with shape $\theta$ and unit scale.  If $X=(X_1,\ldots,X_n)$ consists of iid samples from this model, the relative likelihood function is 
\[ R(X,\theta) = \Bigl\{ \frac{\Gamma(\hat\theta_X)}{\Gamma(\theta)} \Bigr\}^n \Bigl( \prod_{i=1}^n X_i \Bigr)^{\theta - \hat\theta_X}, \]
where $\hat\theta_X$ is the maximum likelihood estimator, given by a solution to the equation 
\[ \psi(\theta) = \frac1n \sum_{i=1}^n \log X_i, \]
with $\psi = (\log \Gamma)'$, the digamma function or derivative of the log gamma function.  There is no closed-form expression for this estimator, but it can be readily evaluated numerically using the {\tt digamma} and {\tt uniroot} functions in R.  Here, we sampled data of size $n=5$ from $\prob_\Theta$, with $\Theta=5$, and a plot of the IM contour function is shown in Figure~\ref{fig:gam}(a).  The distribution function $Q_x^\star$ of the corresponding inner probabilistic approximation is shown in Figure~\ref{fig:gam}(b) along with the Bayesian posterior distribution function based on Jeffreys prior.  For one-parameter problems, the Jeffreys-prior Bayes solution performs quite well, so it should come as no big surprise that the two solutions, despite having entirely different motivations, are virtually the same.  The distribution functions can be differentiated to get the corresponding density functions, and these are shown in Figure~\ref{fig:gam}(c).  While the Bayes posterior density has a mostly closed-form expression (with the exception of the normalizing constant), the density $q_x^\star$ corresponding to $Q_x^\star$ must be evaluated numerically. That evaluation can be be unstable, but one can first approximate $Q_x^\star$ by a suitable spline and then take $q_x^\star$ to be the derivative of that spline approximation, which is what was done to create the plot in Figure~\ref{fig:gam}(c).  
\end{ex}

\begin{figure}[t]
\begin{center}
\subfigure[Possibility contour]{\scalebox{0.39}{\includegraphics{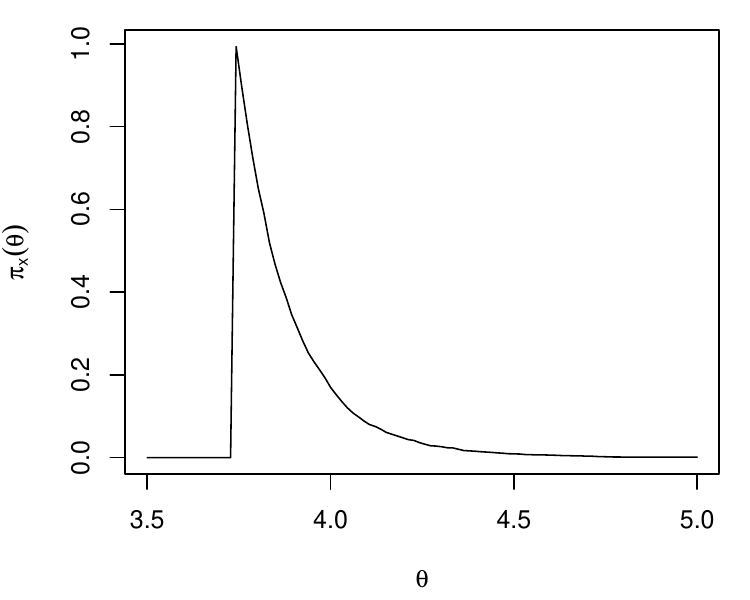}}}
\subfigure[Distribution function]{\scalebox{0.39}{\includegraphics{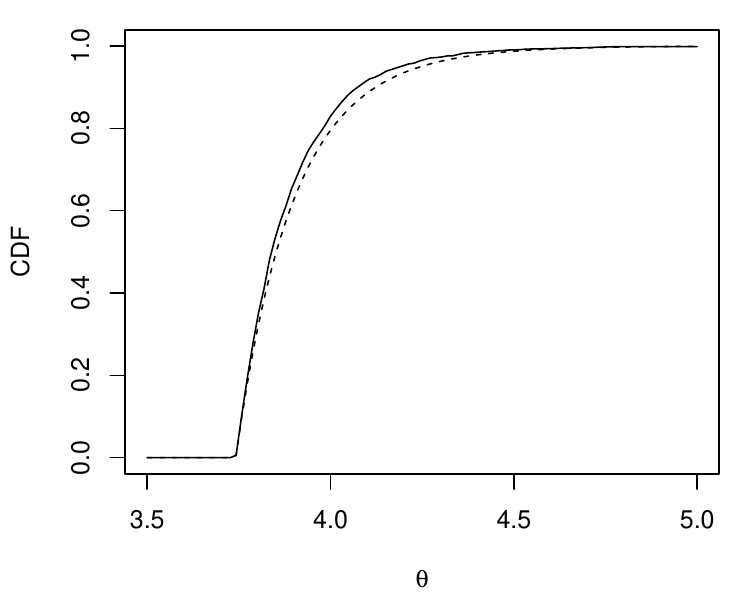}}}
\end{center}
\caption{Results for the uniform illustration in Example~\ref{ex:unif.endpoints}.  Panel (b) compares the IM result (solid) with the Hellinger information-based prior Bayes result (dashed).}
\label{fig:unif.endpoints}
\end{figure}

\begin{ex}
\label{ex:unif.endpoints}
Consider the model $\prob_\theta=\unif(\theta, \theta^2)$, where $\theta > 1$.  This is a classic example of a ``non-regular'' problem in the sense that the distribution's support depends on the unknown parameter.  Therefore, we cannot expect any limiting Gaussianity in this example.  We can, however, construct a possibilistic IM that offers exactly valid uncertainty quantification, just as in the main text.  If $X=(X_1,\ldots,X_n)$ is an iid sample, then the likelihood function is 
\[ L_X(\theta) = (\theta^2 - \theta)^{-n} \, 1(X_{(n)}^{1/2} \leq \theta \leq X_{(1)}), \]
where $X_{(1)}$ and $X_{(n)}$ are the smallest and largest of the sample values, respectively.  Note that, in this case, the pair $(X_{(1)}, X_{(n)})$ is a minimal sufficient statistic for the scalar parameter $\theta$; this difference in dimension is another consequence of the model's non-regularity.  The maximum likelihood estimator is $\hat\theta_X = X_{(n)}^{1/2}$, and the relative likelihood can be readily calculated as 
\[ R(X,\theta) = \left( \frac{X_{(n)} - X_{(n)}^{1/2}}{\theta^2 - \theta} \right)^n \, 1(X_{(n)}^{1/2} \leq \theta \leq X_{(1)}). \]
Both the ranking and validification steps are easy to carry out, yielding the IM solution.  For illustration, we simulate $n=10$ observations from $\prob_\Theta$, with $\Theta=4$.  Figure~\ref{fig:unif.endpoints}(a) shows the contour function, which has a jump discontinuity at the mode, the maximum likelihood estimator $\hat\theta_X = X_{(n)}^{1/2}$, and is strictly decreasing to the right.  In this case, from Corollary~\ref{cor:monotone}, it is easy to get the distribution function of the inner probabilistic approximation---it is just the complement of the plausibility contour to the right of the mode.  A plot of this distribution function is shown in Figure~\ref{fig:unif.endpoints}(b).  For comparison, we show the Bayesian posterior distribution function with prior based on the arguments in \citet{shemyakin2014}; note that, for non-regular problems like this without a Fisher information, defining ``non-informative'' priors is non-trivial.  In any case, as expected, the two distribution functions are similar, but not the same.  Again, an advantage of the IM-driven solution is that there is no prior distribution that needs to be introduced and justified. 
\end{ex}

\begin{figure}[t]
\begin{center}
\subfigure[Possibility contour]{\scalebox{0.39}{\includegraphics{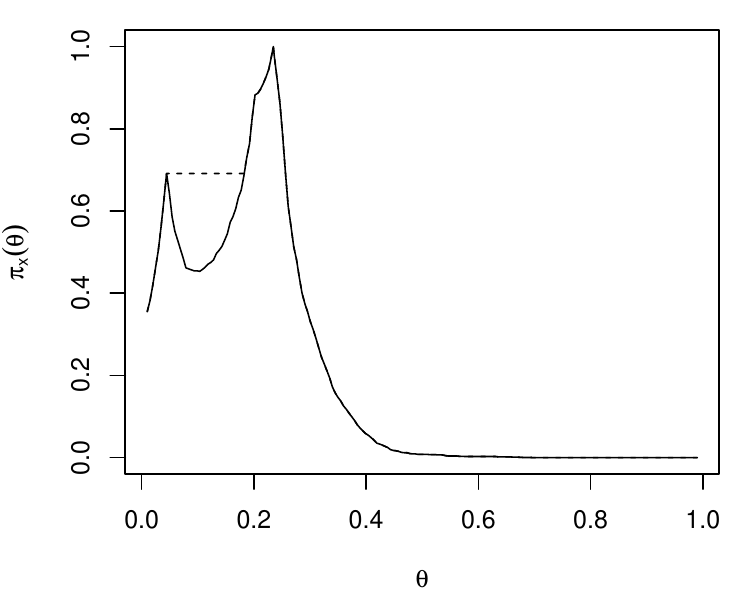}}}
\subfigure[Distribution function]{\scalebox{0.39}{\includegraphics{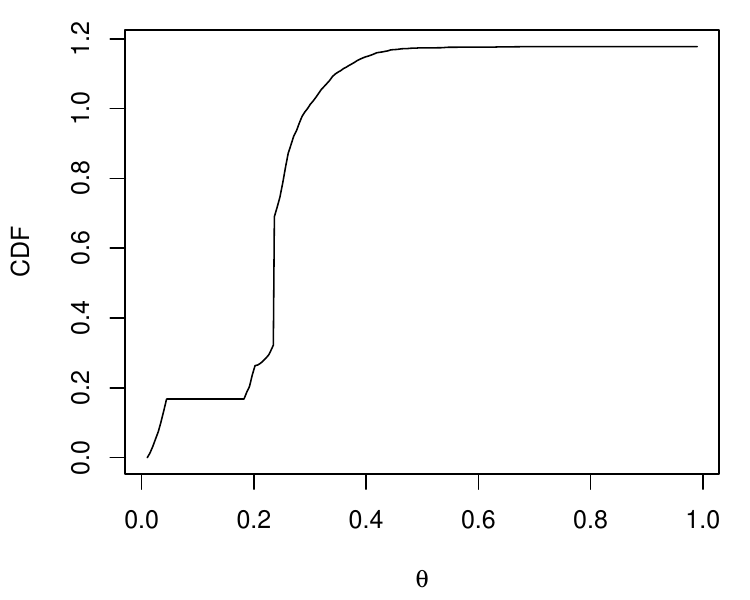}}}
\end{center}
\caption{Results for the triangular model illustration in Example~\ref{ex:tri}.  Panel (a) shows the contour (solid) and its least monotone majorant (dashed).}
\label{fig:tri}
\end{figure}

\begin{ex}
\label{ex:tri}
An interesting non-regular example is the so-called asymmetric triangular distribution \citep[e.g.,][Example~11]{bergerbernardosun2009}, with density function
\[ p_\theta(x) = \begin{cases} 2x/\theta  & \text{if $0 \leq x \leq \theta$}, \\ 2(1-x) / (1-\theta)  & \text{if $\theta < x \leq 1$}, \end{cases} \]
where $\theta \in [0,1]$.  The density has a unique mode at $\theta$, but the graph has a sharp corner and is not differentiable there.  Suppose $X=(X_1,\ldots,X_n)$ is an iid sample and the goal is to make inference on the true $\Theta$.  There are no closed-form expressions that are beneficial to look at here, but suffice it to say that the same possibilistic IM construction applies just as before.  Figure~\ref{fig:tri}(a) shows a plot of the contour based on a simulated sample of size $n=10$ from the model with $\Theta=0.3$.  Two things to notice here are (a) there is an interval of non-monotonicity and (b) the contour does not vanish at the left endpoint.  The non-monotonicity can be easily remedied by working with the least monotone majorant, and this is shown by the dashed line in Figure~\ref{fig:tri}(a).  By replacing the contour by this majorant, we find ourselves in the context covered by Corollary~\ref{cor:boundary} above, and this tells us how to find the inner probabilistic approximation.  The distribution function corresponding to that approximation is shown in Figure~\ref{fig:tri}(b). 
\end{ex}

\bibliographystyle{apalike}
\bibliography{mybib}

\end{document}